\pdfoutput=1
\documentclass[11pt]{article}

\usepackage{amsthm}
\usepackage{amsmath, bm}
\usepackage{amssymb}
\usepackage{natbib}
\usepackage{geometry}
\usepackage{graphicx}
\usepackage{subcaption}
\usepackage{multirow}
\usepackage{listings}
\usepackage{xcolor}
\usepackage{authblk}
\usepackage{float}
\usepackage{subfiles}
\usepackage{indentfirst}
\usepackage{fullpage}

\theoremstyle{definition}
\newtheorem{Def}{Definition}[section]
\newtheorem{Thm}{Theorem}[section]

\lstdefinestyle{CStyle} {
    language=C,
    frame=single,
	backgroundcolor = \color{white},
	basicstyle=\ttfamily,
	keywordstyle=\color{blue},	
    commentstyle=\color{red},
    stringstyle=\color{black},
    directivestyle=\color{purple},
    showspaces=false,
    showstringspaces=false
}

\newcommand{\diag}[1]{{\rm diag}\LRp{#1}}
\newcommand{\td}[2]{\frac{{\rm d}#1}{{\rm d}{\rm #2}}}
\newcommand{\pd}[2]{\frac{\partial#1}{\partial#2}}

\newcommand{\LRp}[1]{\left( #1 \right)}

\newcommand{\LRc}[1]{\left\{ #1 \right\}}

\newcommand{\avg}[1] {\ensuremath{\LRc{\!\{#1\}\!}}}

\title{High-order entropy stable discontinuous Galerkin methods for the shallow water equations: curved triangular meshes and GPU acceleration}
\author[1]{Xinhui Wu}
\author[2]{Ethan J.\ Kubatko}
\author[1]{Jesse Chan}
\affil[1]{Department of Computational and Applied Mathematics, Rice University}
\affil[2]{Department of Civil, Environmental and Geodetic Engineering, The Ohio State University}

\begin{document}
\maketitle

\begin{abstract}
We present a high-order entropy stable discontinuous Galerkin (ESDG) method for the two dimensional shallow water equations (SWE) on curved triangular meshes. The presented scheme preserves a semi-discrete entropy inequality and remains well-balanced for continuous bathymetry profiles. We provide numerical experiments which confirm the high-order accuracy and theoretical properties of the scheme, and compare the presented scheme to an entropy stable scheme based on simplicial summation-by-parts (SBP) finite difference operators. Finally, we report the computational performance of an implementation on Graphics Processing Units (GPUs) and provide comparisons to existing GPU-accelerated implementations of high-order DG methods on quadrilateral meshes. 
\end{abstract}

\section{Introduction}
\label{sec:Intro}
\begin{subequations}
\numberwithin{equation}{section}

The aim of this paper is to present and compare two high-order entropy stable discontinuous Galerkin (DG) schemes. The first is a modal DG formulation of the shallow water equations (SWE), for which the volume and surface quadrature rules can be chosen arbitrarily \cite{chan2018discretely}. The second scheme uses triangular summation-by-parts finite difference operators whose construction is based on carefully chosen quadrature rules which satisfy certain accuracy conditions and contain boundary points \cite{chen2017entropy}. We also implemented both schemes on GPUs for computational acceleration. The analysis and optimization of the multi-threading OCCA code will be discussed. 

The shallow water equations (SWE) are a popular mathematical model for fluid flows in rivers, lakes and coastal regions, where the horizontal scales are much greater than the vertical ones. The shallow water equations in 2D are \cite{wintermeyer2018entropy}
\label{eq:SWE}
\begin{align}
h_t+(hu)_x+(hv)_y&=0,         			 \label{eq:SWEa} \\
(hu)_t+(hu^2+gh^2/2)_x+(huv)_y&=-ghb_x,  \label{eq:SWEb} \\ 
(hv)_t+(huv)_x+(hv^2+gh^2/2)_y&=-ghb_x. \label{eq:SWEc}
\end{align}


The height of the water is denoted by $h = h(x,y,t)$, as measured from the bottom, and should be positive $h > 0$. The velocity in the $x$ direction is denoted by $u = u(x,y,t)$ and the velocity in the $y$ direction is denoted by $v = v(x,y,t)$. The gravitational constant is denoted by $g$. The bathymetry height is denoted by $b=b(x,y)$, and is assumed constant over time. The subscript $(.)_t$ denotes the time derivative, while the subscripts $(.)_x$ or $(.)_y$  denote directional derivatives along $x$ and $y$ axis, respectively. We also define the height of water free surface $H=h+b$. 

The SWE are derived from the Navier-Stokes equations, which describe the motion of incompressible fluids such as water. The Navier-Stokes equations are themselves derived from the equations for conservation of mass and linear momentum. By specifying boundary conditions for the Navier-Stokes equations for a water column and integrating over the depth of the column, one obtains the SWE system \cite{dawson2008shallow}. One of the most prominent practical applications of the SWE is the numerical prediction of storm surges under extreme weather conditions like hurricanes near coastal regions \cite{dawson2011discontinuous, akbar2013hybrid}. 



In addition to being accurate and efficient, numerical methods for the SWE should also preserve certain solutions. Of particular importance is the preservation of steady state solutions, also known as the ``lake at rest'' condition \cite{wintermeyer2018entropy}:
\begin{align}
H=\text{constant},\ u=v=0.
\end{align}
A numerical scheme that preserves this steady state is said to be well-balanced \cite{gassner2016well,noelle2007high}. Schemes that are not well-balanced can generate spurious waves in the presence of varying bottom topographies. A good numerical method for SWE should capture both steady states and their small perturbations so as to avoid the generation of spurious waves \cite{wintermeyer2017entropy}.  Accomplishing this discretely can be challenging, especially for discontinuous bottom topographies, where special discretizations of the source terms are required \cite{wintermeyer2017entropy, fjordholm2011well, xing2014exactly}.

A natural strategy to achieve greater accuracy in numerical simulations is to use higher order methods. This leads to low dissipation errors and long-time accuracy, which are important in simulations of waves. Numerical schemes for the SWE should also address issues unique to the SWE equations, such as well-balancedness and wetting and drying. Finally, since the SWE are non-linear, the solution may become discontinuous even if the initial condition is smooth.  For many numerical methods, stability is a challenge when solving non-linear PDEs, especially in the presence of discontinuous solutions. 

Common numerical methods including the finite difference method, the finite volume method and the finite element method have been applied to the SWE system. This paper focuses on the DG method, which combines advantages of finite element and finite volume methods. DG methods provide a natural path to high-order accuracy and accommodate complex geometries through unstructured meshes. Furthermore, DG methods are simple to parallelize and can take advantage of acceleration using GPU.

Early methods for the SWE were typically low-order accurate and utilized structured grids, which are less geometrically flexible. High-order DG methods address these shortcomings, but introduce issues of stability.  For example, higher order polynomials tend to oscillate in larger magnitude near a discontinuous shock and can result in blow-up of the solution. However, recent work on entropy-stable high-order DG methods \cite{wintermeyer2017entropy, chan2018discretely, wen2020entropy} provide a way to address such instabilities.  Entropy stable DG methods can also be extended to curved meshes, which can be necessary when dealing with complex geometries. 
 
Traditional entropy stable DG formulations have relied on specific finite difference summation-by-parts (SBP) operators, which are constructed using carefully designed quadrature rules which contain boundary points while satisfying certain accuracy conditions. High-order entropy stable DG schemes were more recently extended to ``modal'' formulations, which allow for arbitrary pairings of basis functions and volume/surface quadrature rules. Entropy stable and well-balanced modal DG formulations were introduced for the SWE on Cartesian quadrilateral meshes in \cite{wen2020entropy}. Our goal is to extend these results to curved triangular meshes while ensuring the satisfaction of key properties such as well-balancedness. We examine the numerical behavior of this scheme, and analyze the computational performance of the method with special attention to GPU optimization. 

The new contributions in this paper as follows: in addition to extending entropy stable DG methods for the SWE to curved triangular meshes, we provide a connection between traditional SBP operators and hybridized SBP operators used in ``modal'' entropy stable DG formulations. We also analyze the efficiency of GPU implementations of entropy stable DG methods on triangular meshes. In \cite{wintermeyer2018entropy}, it was shown that GPU implementations of entropy stable methods on quadrilateral meshes do not introduce additional computational cost compared as traditional collocation-type DG schemes. We demonstrate in this work that, while GPU acceleration does provide significant speedups, the cost of entropy stable DG schemes remains higher than the cost of traditional DG schemes on triangular meshes.

This paper begins with a review of entropy stable modal DG schemes for the two dimensional shallow water equations in Section \ref{sec:DG_2dSWE}. In Section \ref{sec:SBPFormulation}, we briefly discuss the SBP formulation and its link with the modal DG formulation. We present numerical results which validate theoretical properties of our schemes in Section \ref{sec:numerical_results}. Section \ref{sec:GPU_opt} provides a description of our GPU implementation and optimization details of the OCCA code. We conclude in Section \ref{sec:Conclusion} with a summary of results.

\end{subequations}

\section{DG method on the 2D shallow water equations}
\label{sec:DG_2dSWE}
\begin{subequations}
\numberwithin{equation}{section}
\subsection{Mathematical assumption and notations}
We first introduce some underlying mathematical assumptions and notations for our DG method. For consistency, we reuse notation from \cite{chan2018discretely}, with slight modifications to provide a cleaner discrete formulation. We denote the triangular reference element by $\hat{D}$ with boundary $\partial{\hat{D}}$. The vertices of the reference triangle are $(-1,-1)$, $(-1,1)$ and $(1,-1)$.  We use $\hat{n}_i$ to represent the $i$th component of the outward normal vector scaled by the face Jacobian on the boundary of the reference element. 
The space of polynomials up to degree $N$ on the reference element is defined as
\begin{align}
P^N(\hat{D}) = \{ \hat{x}^i\hat{y}^j,\quad (\hat{x},\hat{y}) \in \hat{D}, \quad 0 \leq i+j \leq N \}.
\end{align}
Finally, we denote the dimension of the $P^N$ as $N_p = \rm{dim}(P^N(\hat{D}))$.

We wish to build foundations for a discrete matrix-vector formulation of our DG method.  We assume the solution on the reference element $u(\bm{x}) \in P^N(\hat{D})$ such that it can be represented in some polynomial basis $\{\phi_i\}_{i=1}^{N_p}$ of degree up to N, as:
\begin{align}
u(\bm{x}) = \sum_{i=1}^{N_p} u_i \phi_i( \hat{\bm{x}} ), \hspace{1cm} P^N(\hat{D}) = {\rm span} \{\phi_i(\hat{\bm{x}}) \}_{i=1}^{N_p}.
\end{align}
We denote the number of volume and surface quadrature nodes by $N_q$ and $N^f_q$ respectively. Moreover, we assume the volume quadrature rule $\{(\bm{x}_i, w_i )\}_{i=1}^{N_q}$ exactly integrates polynomials of degree at least $(2N-1)$ on the reference element $\hat{D}$ and that the surface quadrature $\{(\bm{x}_i^f, w_i^f )\}_{i=1}^{N_q^f}$ integrates polynomials of degree at least $2N$ on the faces of $\hat{D}$.

Let $\bm{W}$ denote the diagonal $N_q \times N_q$ matrix whose entries are $\bm{W}_{ii} = w_i$, where $w_i > 0$ corresponds to volume quadrature weight. We also define the diagonal matrix $\bm{W}_{f} = {\rm diag}(w_i^f)$ for the surface quadrature weights. We then define the volume and surface quadrature interpolation matrices $\bm{V}_q$ and $\bm{V}_f$ as:
\begin{align}
(\bm{V}_q)_{ij} &= \phi_j(\hat{\bm{x}}_i),\quad 1 \leq j \leq N_p,\quad  1 \leq i \leq N_q, \\
(\bm{V}_f)_{ij} &= \phi_j(\hat{\bm{x}}_i^f),\quad 1 \leq j \leq N_p,\quad  1 \leq i \leq N_q^f,
\end{align}

Let $f(\bm{x})$ denote some polynomial in the basis $\phi_j$ with coefficients $\bm{f}$. The matrix $\bm{V}_q$ maps coefficients $\bm{f}$ to evaluations of $f(\bm{x})$ at volume quadrature points and, similarly, the matrix $\bm{V}_f$ interpolates ${f}$ to surface quadrature points. For example:
\begin{align}
\bm{f}_q = \bm{V}_q\bm{f}, \qquad (\bm{f}_q)_i = f(\hat{\bm{x}}_i),\qquad 1\leq i \leq N_q.
\end{align}

We now define $\bm{D}^i$ as the differentiation matrix with respect to the $i$th coordinate. We may denote the matrices $\bm{D}^1, \bm{D}^2$ as  $\bm{D}^x$ and $\bm{D}^y$ in two-dimensional case. $\bm{D}^i$ is defined implicitly with:
\begin{align}
u(\bm{x}) = \sum_{i=1}^{N_p} u_i \phi_i( \hat{\bm{x}} ), \hspace{1cm} \frac{\partial u}{\partial \hat{x}_i} = \sum_{j=1}^{N_p} (\bm{D}^i\bm{u})_j \phi_j( \hat{\bm{x}} ).
\end{align}
$\bm{D}^i$ maps the basis coefficients of some polynomial $\bm{u} \in P^N$ to coefficients of its $i$th directional derivative with respect to the reference coordinate $\bm{x}_i$.

With the matrix $\bm{V}_q$, we can now introduce the element mass matrix whose entries are the evaluations of inner products of different basis functions with quadrature points:
\begin{align}
\bm{M} = \bm{V}_q^T\bm{W}\bm{V}_q, \hspace{0.5cm} \bm{M}_{ij} = \sum_{k=1}^{N_q} w_k\phi_j(\bm{\hat{x}}_k)\phi_i(\bm{\hat{x}}_k) \approx \int_{\hat{D}} \phi_j\phi_i d\bm{\hat{x}} = (\phi_j,\phi_i)_{\hat{D}}.
\end{align}

We also define a $L^2$ projection operator $\Pi_N: L^2(\hat{D}) \to P^N(\hat{D})$ such that
\begin{align}
(\Pi_N f, v)_{\hat{D}} = (f,v)_{\hat{D}},\hspace{1cm} \forall v \in P^N(\hat{D}).
\end{align}
When integrals within the $L^2$ projection are computed with quadrature, the discrete quadrature-based $L^2$ projection of a function $f(x)$ can be expressed as following:
\begin{align}
\bm{M}\bm{u} = \bm{V}^T_q\bm{W}\bm{f}, \hspace{1cm} \bm{f}_i = f(\hat{\bm{x}}_i),\hspace{1cm} 1\leq i \leq N_q,
\end{align}
where $\bm{u}$ is the vector of coefficients of the quadrature-based $L^2$ projection of the function values of $\bm{f}$. We can define the quadrature-based $L^2$ projection matrix $\bm{P}_q$, by inverting the mass matrix:
\begin{align}
\bm{P}_q = \bm{M}^{-1}\bm{V}_q^T\bm{W}.
\end{align}
The matrix $\bm{P}_q$ maps a function in terms of its evaluations at quadrature points to its coefficients of the $L^2$ projection in the basis $\phi_i(\bm{\hat{x}})$.
Notice that since $\bm{M} = \bm{V}_q^T\bm{W}\bm{V}_q$, we have
\begin{align}
\bm{P}_q\bm{V}_q = \bm{M}^{-1}\bm{V}_q^T\bm{W}\bm{V}_q = \bm{I}. \label{eq:PqVq}
\end{align}
This implies that when we apply $\bm{P}_q$ to the evaluations of polynomial function at volume quadrature points, we recover the coefficients of the polynomial in the basis $\phi_i(\bm{\hat{x}})$.

\subsection{Discrete formulation in 2D}
With the tools defined in the previous sections, we can now derive discretization matrices which will be used in our discrete DG formulation.

In d-dimensions, define the following matrices
\begin{align}
\bm{\hat{Q}}^i = \bm{M}\bm{D}^i, \hspace{1cm} \bm{B}^i = \bm{W}_f \diag{\bm{\hat{n}}_i}, \hspace{1cm} i=1,...,d.
\end{align}
With the above definitions, we have that
\begin{align}
\label{eq:SBP_Q}
\bm{\hat{Q}}^i + (\bm{\hat{Q}}^i)^T = \bm{V}_f^T\bm{B}^i\bm{V}_f.
\end{align}

By combining the projection matrix $\bm{P}_q$ with the matrix $\bm{\hat{Q}}^i$, we can construct a nodal differentiation operator at quadrature points \cite{chan2018discretely}:
\begin{align}
\bm{Q}^i = \bm{P}_q^T\bm{\hat{Q}}^i\bm{P}_q.
\end{align}

We also define the the matrix $\bm{E}$, which extrapolates volume quadrature nodes to surface quadrature nodes, as 
\begin{align}
\bm{E} = \bm{V}_f\bm{P}_q.
\end{align}
Then we have the following generalized SBP property:
\begin{align}
\bm{Q}^i + (\bm{Q}^i)^T = \bm{E}^T\bm{B}^i\bm{E}.  
\end{align}
Similarly, for convenience, we define $\bm{V}_h$ as
\begin{align}
\bm{V}_h =
\begin{bmatrix}
\bm{V}_q\\
\bm{V}_f\\
\end{bmatrix}.
\end{align}

\subsection{Hybridized SBP operators}\label{subsec:Hybridized_SBP_operators}
Entropy stable formulations for nonlinear conservation laws can be constructed using the generalized SBP operators introduced in the previous section. However, these schemes introduce numerical flux terms which couple together all degrees of freedom on neighboring elements \cite{crean2017high}. 

To avoid this, we introduce the hybridized operator $\bm{Q}_h^i$, which is given explicitly as 
\begin{align}
\bm{Q}_h^i = \frac{1}{2}
\begin{bmatrix}
\bm{Q}^i - (\bm{Q}^i)^T & \bm{E}^T\bm{B}^i\\
-\bm{B}^i\bm{E} & \bm{B}^i\\
\end{bmatrix}.
\end{align}
This operator is designed to be applied to vectors of solution values at both volume and surface quadrature nodes and mimics the structure of boundary terms used in hybridized DG methods \cite{chen2019review}.  When used in a DG formulation, it allows one to construct entropy stable formulations using more standard DG numerical fluxes. 

We have the following theorem:
\begin{Thm}
$\bm{Q}_h^i$ satisfies the $SBP-like$ property \cite{chan2018discretely}:
\begin{align}
\bm{Q}_h^i +\LRp{\bm{Q}_h^i}^T = \bm{B}_h^i, \hspace{1cm} \bm{B}_h^i = \begin{bmatrix}
\bm{0} &  \\
  &  \bm{B}^i\\
\end{bmatrix},
\end{align}
 and $\bm{Q}_h^i \bm{1} = 0$, where $\bm{1}$ is the vector of all ones.
\label{thm:SBP_Q_h^i}
\end{Thm} 

\subsection{SWE entropy and entropy variables}
 In this section, we introduce the entropy function and associated entropy variables for the SWE. 
Entropy stability \cite{tadmor2003entropy} is the extension of $L^2$ (energy) stability for linear hyperbolic PDEs to nonlinear conservation laws.  To provide a statement of entropy stability, we first need to define a convex entropy function $S(\bm{u})$. Solutions to nonlinear conservation laws are typically non-unique.  To determine unique solutions, we require that solutions satisfy an entropy inequality. 

The entropy function for the SWE is the total energy of the system \cite{fjordholm2011well, wintermeyer2017entropy}:
\begin{align}
S(\bm{u}) = \frac{1}{2}h(u^2+v^2) + \frac{1}{2}gh^2+ghb.
\label{eq:entropy_function}
\end{align}
We also define the entropy variable $\bm{v} = S'(\bm{u})$. The convexity of the entropy function guarantees that the mapping between $\bm{u}$ and $\bm{v}$ is invertible.  The entropy variables for the SWE are given explicitly as:
\begin{align}
v_1 &= \frac{\partial S}{\partial h}  = g(h+b) - \frac{1}{2}u^2 - \frac{1}{2}v^2,\\
v_2 &= \frac{\partial S}{\partial (hu)} = u,\\
v_3 &= \frac{\partial S}{\partial (hv)} = v.
\end{align}
It can be shown as in \cite{mock1980systems} that there exists an entropy flux function $F(\bm{u})$ and entropy potential $\psi(\bm{u})$ such that
\begin{align}
\bm{v}(\bm{u})^T\frac{\partial \bm{f}}{\partial \bm{u}} = \frac{\partial F(\bm{u})^T}{\partial \bm{u}}, \hspace{.5cm}
\psi(\bm{u}) = \bm{v}(\bm{u})^T\bm{f}(\bm{u}) - F(\bm{u}), \hspace{0.5cm} \psi'(\bm{u}) = \bm{f}(\bm{u}).
\end{align}
Assuming for simplicity that bathymetry is constant and that the domain $\Omega$ is periodic, an entropy equality can be derived for smooth solutions $\bm{u}$ by multiplying the SWE by $\bm{v}^T$ and integrating over the domain. Then, using the chain rule and definition of the entropy flux, we have the following statement of entropy conservation
\begin{align}
\int_{\Omega}\frac{\partial S(\bm{u})}{\partial t} = \bm{0}.
\label{eq:ec_bottom0}
\end{align} 
For viscosity solutions, it can be shown that (\ref{eq:ec_bottom0}) becomes an entropy inequality.

Our goal is to reproduce this statement of entropy conservation discretely. The resulting entropy conservative formulation can then be used to construct entropy stable formulations by adding appropriate entropy dissipation terms.

\subsection{Entropy conservation and  flux differencing}\label{sec:EC_flux}
In this section, we introduce numerical fluxes for SWE  and describe an entropy conservation discrete formulation \cite{fisher2013high, gassner2018br1, chen2017entropy, gassner2016split}. To construct the entropy stable scheme in higher dimensions, we require entropy conservative fluxes as defined in \cite{tadmor1987numerical}:

\begin{Def}{}
Let $\bm{f}_S^i(\bm{u}_L, \bm{u}_R)$ be a bivariate function which is symmetric and consistent with the flux function $\bm{f}^i(\bm{u})$, for $i = 1,...,d$
\begin{align}
\bm{f}_S^i(\bm{u},\bm{u}) = \bm{f}^i(\bm{u}), \hspace{1cm} \bm{f}_S^i(\bm{u},\bm{v}) = \bm{f}_S^i(\bm{v},\bm{u}).
\end{align}
The numerical flux $\bm{f}_S^i(\bm{u}_L, \bm{u}_R)$ is entropy conservative if, for entropy variable $\bm{v}_L = \bm{v}(\bm{u}_L)$, $\bm{v}_R = \bm{v}(\bm{u}_R)$
\begin{align}
\left( \bm{v}_L - \bm{v}_R\right)^T \bm{f}_S^i\left(\bm{u}_L,\bm{u}_R\right) = \psi^i_L - \psi^i_R, \\
\psi^i_L = \psi^i(\bm{v}(\bm{u}_L)), \quad \psi^i_R = \psi^i(\bm{v}(\bm{u}_R)).
\end{align}
\label{def:consevative_flux}
\end{Def}
The flux $\bm{f}_S^i$ can be used to construct entropy conservative and entropy stable finite volume methods. Entropy stable finite volume schemes were generalized in \cite{fjordholm2012arbitrarily} to arbitrary high order. 

This numerical flux is also used for the construction of discretely entropy stable DG schemes using an approach referred to as flux differencing \cite{fisher2013high, carpenter2014entropy,gassner2018br1, chen2017entropy}. Flux differencing was first used to systematically recover entropy stable split formulations in \cite{gassner2016split}, but is applicable to a broader range entropy stable formulations.  Entropy stable DG schemes also couple elements together using using the same entropy conservative flux $\bm{f}_S^i(\bm{u}_L, \bm{u}_R)$ as an interface flux \cite{carpenter2014entropy, gassner2018br1, chen2017entropy}. Entropy stable schemes are typically constructed by first constructing an entropy conservative scheme, then adding entropy dissipation through appropriate penalization terms at element interfaces. These additional penalization terms convert schemes which satisfy a global entropy equality into schemes which satisfy a global entropy inequality. 


Using flux differencing from \cite{chan2018discretely, crean2018entropy}, we can replace the term $\bm{f}^i(\bm{u}(x))$ with the term $2\bm{f}_S^i(\bm{u}(x), \bm{u}(x))$. Then we define a flux matrix $\bm{F}^i$ as the evaluations of $\bm{f}_S^i(\bm{u}(x), \bm{u}(y))$ at quadrature points:
\begin{align}
(\bm{F}^i)_{jk} = \bm{f}_S^i(\bm{u}(\hat{x}_j), \bm{u}(\hat{x}_k)), \qquad 1 \leq j,k \leq N_q.
\label{eq:flux_matrix}
\end{align}
The term $2(\bm{Q}\circ\bm{F})\bm{1}$ approximates $\int \frac{\partial \bm{f}^i(\bm{u}(x))}{\partial x}$, and the key idea in entropy stable DG formulations is to replace $\bm{Q}\bm{f}(\bm{u})$ with $2(\bm{Q}\circ\bm{F})\bm{1}$, where $\bm{Q}\circ\bm{F}$ denotes the Hadamard product between $\bm{Q}$ and $\bm{F}$.

\subsection{Entropy projection}
We seek a degree $N$ polynomial approximation of the conservative variables $\bm{u}(x,t)$ with coefficients $\bm{u}_h(t)$ such that
\begin{align}
\bm{u}_N(\hat{\bm{x}},t) = \sum_{i=1}^{N_p} (\bm{u}_h(t))_i \phi_i(\hat{\bm{x}}), \hspace{1cm} (\bm{u}_h(t))_i \in \mathbb{R}^n.
\end{align}
Because $\bm{u}_h$ consists of vectors of coefficients for each scalar component of $\bm{u}_N(\hat{\bm{x}},t)$, we should understand the discretization matrices as being applied to vectors like $\bm{u}_h$ in a Kronecker product sense. For example, $\bm{A}\bm{u}_h$ should be interpreted as applying $\bm{A}$ to each component of $\bm{u}_h$.

Reproducing conservation of entropy discretely faces an additional challenge.  Since the SWE system is non-linear, entropy variables are not contained in the approximation space, and in general, $\bm{v}(\bm{u}) \not\in P^N$ even if $\bm{u}\in P^N$.  For shallow water, if one assumes $ h, hu \in P^N$, then  $v_2(h,hu) = u = \frac{hu}{h}$ is rational and non-polynomial. 

Unfortunately, for DG methods, the test space contains only piecewise polynomial functions.  To circumvent this issue, we introduce $\bm{v}_h$ as the $L^2$ projection of the entropy variables and $\bm{\tilde{u}}$ as the evaluations of the conservative variables in terms of the $L^2$ projected entropy variables
\begin{align}
\bm{u}_q = \bm{V}_q\bm{u}_h, \hspace{1cm} \bm{v}_q = \bm{v}(\bm{u}_q), \hspace{1cm} \bm{v}_h = \bm{P}_q\bm{v}_q, 
\end{align}
\begin{align}
\bm{\tilde{v}} = 
\begin{bmatrix}
\bm{\tilde{v}}_q\\ \bm{\tilde{v}}_f
\end{bmatrix}
 = \begin{bmatrix}
 \bm{V}_q\\ \bm{V}_f
 \end{bmatrix}
 \bm{v}_h, \hspace{1cm}
 \bm{\tilde{u}} = 
\begin{bmatrix}
\bm{\tilde{u}}_q\\ \bm{\tilde{u}}_f
\end{bmatrix} = 
 \bm{u}(\bm{\tilde{v}}).
\end{align}
Here $\bm{u}_q$ and $\bm{v}_q$ denote the conservative variables and entropy variables evaluated at the volume quadrature points. The vector $\bm{\tilde{v}}$ denotes the evaluations of the $L^2$ projection of the entropy variables at both volume and surface quadrature points, while $\bm{\tilde{u}}$ denotes the evaluations of the conservative variables in terms of the projected entropy variables $\bm{u}(\Pi_N \bm{v})$.

\subsection{Curved triangular meshes}\label{sec:curved_mesh}
We now extend the construction of hybridized SBP operators to curved triangular meshes, where each element can be represented as a curvilinear mapping $\Phi^k$ of the reference element $\hat{D}$ \cite{chan2018discretely, ranocha2017extended}. Because we map the reference triangle to curved triangular elements, the geometric factors are not constant anymore, which can impact the stability of our DG formulation.  We want to ensure the entropy stability on curved triangular meshes. 

We introduce the following definitions associated with Jacobians of the mapping $\Phi^k$:
\begin{itemize}
\item$J^k$ denotes the determinant of the Jacobian of $\Phi^k$.
\item$\bm{J}^k_f$ denotes the vector of $J^k_f$ at surface quadrature points.
\item$\bm{\hat{J}}_f$ denotes the vector contains $\hat{J}_f$, the face Jocabian factor of the mapping from faces of the reference elements to the reference face. We assume $\hat{J}_f$ is pre-multiplied into the surface quadrature weights.
\end{itemize}
We introduce the ``split'' form of derivative to preserve entropy stability\cite{gassner2016split, chan2018discretely, ranocha2017extended}
\begin{align}
\pd{u}{x_i} = \frac{1}{2} \sum_j \LRp{\pd{\hat{x}_j}{x_i}\pd{u}{\hat{x}_j} + \pd{}{\hat{x}_j}\LRp{u\pd{\hat{x}_j}{x_i}}}.
\end{align}
Let $\bm{G}^k_{ij} = \frac{\partial \hat{x}_j}{\partial x_i}$ be the vector of geometric factors evaluated at the quadrature points on element $D^k$. We define the physical SBP operator on $D^k$ as:
\begin{align}
\bm{Q}_h^{i,k} = \frac{1}{2}\sum_{j=1}^{2} \diag{\bm{G}^k_{ij}}\bm{\hat{Q}}_h^{j,k}+\bm{\hat{Q}}_h^{j,k}\diag{\bm{G}^k_{i,j}},
\end{align}
It can be shown that $\bm{Q}^{i,k}_h$ satisfies the following SBP property on $D^k$
\[
\bm{Q}^{i,k}_h + \LRp{\bm{Q}^{i,k}_h}^T = \begin{bmatrix}
\bm{0} &\\
&\bm{B}^{i,k}\end{bmatrix}, \qquad \bm{B}^{i,k} = \bm{W}_f \diag{\bm{J}_f^k/\bm{\hat{J}}_f \circ \bm{n}^k_i}, 
\]
where  $\bm{J}_f^k/\bm{\hat{J}}_f$ denotes element-wise division between the vectors containing the face Jacobians and the reference face Jacobians, and $\bm{n}^k_i$ is a vector of the outward unit normals on $D^k$. 

Finally, we define the mass matrix for the element $D^k$
\begin{align}
\bm{M}_h^k = \bm{V}_q^T\bm{W}^k\bm{V}_q = \bm{V}_q^T\bm{W}J^k\bm{V}_q,
\end{align}
where the $L^2$ mass matrix is now weighted by a non-constant Jacobian $J^k$.

\subsection{Entropy conservative flux for SWE}

In this section, we finalize our discrete DG formulation for the SWE. We first present entropy conservative (EC) fluxes for the 2D SWE \cite{tadmor2003entropy, carpenter2014entropy, gassner2016well, wintermeyer2018entropy}
\begin{align}
\bm{f}^x_{S}\LRp{\bm{u}_L,\bm{u}_R} &=
\begin{bmatrix}
\avg{hu}\\
\avg{hu}\avg{u} + g\avg{h}^2 - \frac{1}{2}g\avg{h^2}\\
\avg{hu}\avg{v}
\end{bmatrix}
\label{eq:SWE_flux1}
, \\
\bm{f}^y_{S}\LRp{\bm{u}_L,\bm{u}_R} &=
\begin{bmatrix}
\avg{hv}\\
\avg{hv}\avg{u}\\
\avg{hv}\avg{v} + g\avg{h}^2 - \frac{1}{2}g\avg{h^2}\\ 
\end{bmatrix}.
\label{eq:SWE_flux2}
\end{align}
A discrete entropy conservative formulation of the SWE is then given as follows on an element $D^k$:
\begin{gather}
\bm{M}_h^k\td{\bm{u}}{t} + \sum_{i=x,y} \begin{bmatrix}
\bm{V}_q\\\bm{V}_f
\end{bmatrix}^T\LRp{2\bm{Q}^{i,k}_{h}\circ \bm{F}^i}\bm{1} + \bm{V}_f^T\bm{B}^{i,k}\LRp{\bm{f}^{i}_S(\tilde{\bm{u}}^+,\tilde{\bm{u}})-\bm{f}^i(\tilde{\bm{u}}_f)} = \bm{S}, \label{eq:dgform}\\
\nonumber
(\bm{F}^i)_{j,k} = \bm{f}^{i}_S(\tilde{\bm{u}}_i,\tilde{\bm{u}}_j),\qquad 1\leq j,k\leq N_q+N_q^f,
\end{gather}
where $\bm{S}$ is the source term
\[
\bm{S} = -g\begin{bmatrix}
\bm{0}\\
{\rm diag}\LRp{\bm{h}}\bm{Q}^x\bm{b}\\
{\rm diag}\LRp{\bm{h}}\bm{Q}^y\bm{b}
\end{bmatrix}.
\]
Recall that $\tilde{\bm{u}}$ are the ``entropy-projected'' conservative variables introduced in the previous section. This formulation was shown to be entropy conservative in \cite{chan2019discretely,chan2019skew}.  

We note that, in our implementation, we precompute the inverses of the element mass matrices, i.e.\ $(\bm{M}_h^k)^{-1}$ for all $k$, and store them on the GPU. This can be avoid using weight-adjusted mass matrix inverses \cite{chan2019discretely}, which we will investigate in future work.


Now we present a proof of entropy conservation for our discrete DG formulation of the SWE.  A proof of entropy conservation for general nonlinear conservation laws was given in \cite{chan2018discretely}, but did not account for bathymetric source terms. This proof extends this theory while accounting for the presence of varying bathymetry. 
\begin{Thm}
Let $\bm{f}_s$ be an entropy conservative flux from Definition \ref{def:consevative_flux}. Then assuming continuity in time, the semi-discrete formulation (\ref{eq:dgform})
is entropy conservative and well-balanced for $b \in P^N$. 
\label{thm:entropy_conservative}
\end{Thm}
\begin{proof}
First, we divide the entropy variable $\bm{v}$ into two parts
\begin{align}
\bm{v} = \bm{v}_0+\bm{v}_b , \hspace{0.5cm} 
\bm{v}_0 = \begin{bmatrix}
gh-\frac{1}{2}(u^2+v^2)\\
u\\
v
\end{bmatrix},\hspace{0.5cm} 
\bm{v}_b = \begin{bmatrix}
gb\\
0\\
0
\end{bmatrix},
\end{align}
where $\bm{v}_0$ are terms in the entropy variables corresponding to zero bathymetry and $\bm{v}_b$ are terms corresponding to the contribution from variable bathymetry. 

Without loss of generality, we assume the solution and $b$ are constant along the $y$ direction such that only derivatives in the $x$-direction remain. The general proof can be treated by repeating this procedure in each coordinate direction. 

Without $y$-derivatives, the formulation (\ref{eq:dgform}) reduces to 1D formulation
\[
\bm{M}\td{\bm{u}}{t} + 2\LRp{\bm{Q}^{x,k}_h\circ \bm{F}_S^x}\bm{1} + \bm{B}^{x,k}\LRp{\bm{f}^x_S\LRp{\tilde{\bm{u}}^+_f,\tilde{\bm{u}}_f}-\bm{f}^x(\tilde{\bm{u}}_f)} = \bm{0}.
\]
From here, we drop superscripts, subscripts, tildes, and $(\cdot)^{x,k}_h$ on all solution variables to simplify notation. 

Notice that the term $2\LRp{\bm{Q}\circ \bm{F}_S}\bm{1}$ is independent of the bathymetry $b$. For $b=0$, multiplying the discrete formulation (\ref{eq:dgform}) by $\bm{v}_0^T$ yields that
\begin{align}
\bm{v}_0^T \LRp{\bm{M}\td{\bm{u}}{t} + 2\LRp{\bm{Q}\circ \bm{F}_S}\bm{1} + \bm{B}\LRp{\bm{f}_S-\bm{f}(\bm{u})}} = \bm{v}_0^T\bm{M}\td{\bm{u}}{t} = 0.
\end{align}
The proof is the same as the one given in \cite{chan2018discretely} (Theorem 2).  

It remains to show that entropy conservation still holds if $b$ is not constant. Multiplying (\ref{eq:dgform}) by $\bm{v}^T$ then yields
\begin{align}
\nonumber
&\bm{v}^T\LRp{\bm{M}\td{\bm{u}}{t} + 2\LRp{\bm{Q}\circ \bm{F}_S}\bm{1} + \bm{B}\LRp{\bm{f}_S-\bm{f}(\bm{u})}} \\
= &\bm{v}_0^T\bm{M}\td{\bm{u}}{t} + \bm{v}_b^T\LRp{\bm{M}\td{\bm{u}}{t} + 2\LRp{\bm{Q}\circ \bm{F}_S}\bm{1} + \bm{B}\LRp{\bm{f}_S-\bm{f}(\bm{u})}} =  \bm{v}^T\bm{s}.
\label{eq:wb_pf_separate_terms}
\end{align}
We wish to show that this implies 
\[
\bm{1}^T\bm{W}\td{S(\bm{u})}{t} = 0.
\]
The time derivative terms $\bm{v}_b^T\bm{M}\td{\bm{u}}{t}$ and $\bm{v}_0^T\bm{M}\td{\bm{u}}{t}$ from (\ref{eq:wb_pf_separate_terms}) combine to the form the term $\bm{v}^T\bm{M}\td{\bm{u}}{t} = \bm{1}^T\bm{W}\frac{\partial S}{\partial t}$ for varying bathymetry (the proof can be found in \cite{chan2018discretely}, and utilizes properties of the $L^2$ projection). What remains to show is then that
\[
\bm{v}_b^T\LRp{2\LRp{\bm{Q}\circ \bm{F}_S}\bm{1} + \bm{B}\LRp{\bm{f}_S-\bm{f}(\bm{u})}} -  \bm{v}^T\bm{s} = 0.
\]
Since the second and the third components of $\bm{v}_b$ are 0, the first expression corresponds to testing the mass conservation equation in the DG formulation with $gb$. Substituting the flux from (\ref{eq:SWE_flux1}), this yields
\begin{align}
&\bm{v}_b^T 2 \LRp{\bm{Q}\circ \bm{F}_S}\bm{1} + \bm{v}_b^T \bm{B} \LRp{\bm{f}_S-\bm{f}(\bm{u})} - (\bm{v}_b^T + \bm{v}_0^T)\bm{s}\\
\nonumber= &g\bm{b}^T\bm{Q}(\bm{hu}) + g(\bm{hu})^T\bm{Q}\bm{b} + \frac{1}{2}g\bm{b}^T\bm{B}[[\bm{hu}]].
\end{align}
where we have expanded out the interface term $\bm{v}_b^T\bm{B}\LRp{\bm{f}_S - \bm{f}}$ to yield an expression involving the jump $[[.]]$ of a quantity across the element boundary.

Using the SBP property yields
\begin{align}
\nonumber&g\bm{b}^T\bm{Q}(\bm{hu}) + g(\bm{hu})^T\bm{Q}\bm{b} + \frac{1}{2}g\bm{b}^T\bm{B}[[\bm{hu}]]\\
\nonumber= &g\bm{b}^T\bm{Q}(\bm{hu}) + \frac{1}{2}g\bm{b}^T\bm{B}[[\bm{hu}]] + g\bm{b}^T(\bm{Q})^T(\bm{hu})\\
= &g\bm{b}^T\bm{B}(\bm{hu}) + \frac{1}{2}g\bm{b}^T\bm{B}[[\bm{hu}]].
\label{eq:ec_boundary_terms}
\end{align}
Consider two neighboring elements $D^k$ and $D^{k+}$. Then, on each element, the boundary terms are
\begin{align}
\nonumber D^k :  &\hspace{0.5cm}g\bm{b}^T\bm{B}^x(\bm{hu}) + \frac{1}{2}g\bm{b}^T\bm{B}^x((\bm{hu})^+ - (\bm{hu})),\\
\nonumber D^{k+}: &\hspace{0.5cm}g\bm{b}^T\bm{B}^{x+}(\bm{hu})^+ + \frac{1}{2}g\bm{b}^T\bm{B}^{x+}((\bm{hu}) - (\bm{hu})^+)\\
&= -g\bm{b}^T\bm{B}^{x}(\bm{hu})^+ - \frac{1}{2}g\bm{b}^T\bm{B}^{x}((\bm{hu})^+ - (\bm{hu})).
\end{align}
Since outward normals are equal and opposite across an interface, $B^{x,+} = -B^{x}$ and the boundary terms (\ref{eq:ec_boundary_terms}) cancel out when summed up over two neighboring elements. This implies that
\begin{align}
\sum_{\partial D^k} g\bm{b}^T\bm{B}^x(\bm{hu}) + \frac{1}{2}g\bm{b}^T\bm{B}^x[[\bm{hu}]] = 0.
\end{align}
Combining the results from the $b=0$ part, we obtain the entropy conservative property.

To show this formulation is well-balanced, we consider the steady state solution
\begin{align}
    h(x,t) = c - b(x), \qquad u(x,t) = 0,
\end{align}
where $c$ is a constant total water height and $b$ is continuous. We immediately obtain well-balancedness for the first equation involving $h$ because $hu = 0$. To see $u(x,t) = 0$ is preserved, we substitute the nonzero terms of flux from (\ref{eq:SWE_flux2}) into formulation (\ref{eq:dgform})
\begin{align}
&2g\LRp{\bm{Q}\circ\LRp{ \avg{\bm{h}}^2 -\frac{\avg{\bm{h}^2}}{2}} }\bm{1} +g\bm{B}\LRp{\avg{\bm{h}}^2 -\frac{\avg{\bm{h}^2}}{2} -\frac{\bm{h}^2}{2}}\\
= & 2g\LRp{\bm{Q}\circ\frac{\bm{h}_i\bm{h}_j}{2}}\bm{1}+g\bm{B}\bm{0}\\
= & g\LRp{\bm{Q}\circ\LRp{\bm{h}_i\bm{h}_j}}\bm{1} \\
= & g\sum_j \bm{Q}_{ij}\bm{h}_i\bm{h}_j, \qquad i = 1,..., N_q. 
\end{align}
Since $\bm{Q}$ differentiate polynomial exactly up to degree $N$, $\bm{Q}\bm{1} = \bm{0}$ and $\bm{Q}\bm{c} = \bm{0}$ for any constant vector $\bm{c}$. We also have the source term
\begin{align}
&-g\cdot{\rm diag}\LRp{\bm{h}}\bm{Q}\bm{b}\\
= &-g\cdot{\rm diag}\LRp{\bm{h}}\bm{Q}\LRp{\bm{c} - \bm{h}}\\
= &g\cdot{\rm diag}\LRp{\bm{h}}\bm{Q}\bm{h}\\
= &g\bm{h}_i \sum_j \bm{Q}_{ij}\bm{h}_j\\
= &g\sum_j \bm{Q}_{ij}\bm{h}_i\bm{h}_j,  \qquad i = 1,..., N_q.
\end{align}
Therefore the volume and surface contributions cancel out with the source term to achieve well-balancedness.
\end{proof}
\color{black}
Thus, our DG formulation is entropy conservative. The fact that the entropy function for the SWE is the total energy of the system also provides a connection between the mathematical stability of our numerical scheme and physical principles. 

To construct an entropy stable scheme, we add entropy dissipative interface penalization terms to the entropy conservative formulation. In this work, we utilize local Lax-Friedrichs penalization \cite{wintermeyer2017entropy, chen2017entropy} applied to the entropy-projected conservative variables. We note that this choice of penalization also preserves the well-balanced property: for a lake-at-rest condition with zero velocity, the entropy projected conservative variables are the same as the original conservative variables. Lax-Friedrichs penalization adds a scaling of the jumps of the conservative variables, and since the jumps of the conservative variables vanish for continuous water height and bathymetry, the entropy stable scheme reduces to the well-balanced entropy conservative scheme for the lake-at-rest condition.

Finally, while we have restricted ourselves to continuous bathymetry for this paper, it is possible to construct well-balanced and entropy stable schemes for discontinuous bathymetry by adding additional interface terms \cite{fjordholm2011well, wintermeyer2017entropy}.

\subsection{Imposing reflective (wall) boundary conditions}

Proofs of entropy stability have involved periodic boundary conditions.  However, one can also show entropy conservation or entropy stability under reflective wall boundary conditions.  To impose reflective ``wall" boundary conditions in an entropy conservative or entropy stable fashion, we follow procedures outlined in \cite{wintermeyer2017entropy, chen2017entropy}.  For surfaces that correspond to the wall, we set $u^+_{\bm{n}} = -u_{\bm{n}}$ and $h^+=h$, where $(.)^+$ denotes the exterior value of the solution used in a numerical flux (e.g., the value of the solution on a neighboring element for an interior interface), and the subscript $(.)_{\bm{n}}$ denotes the normal component of the vector with respect to the wall. 

\end{subequations}
\section{Entropy stable formulations on curved triangular meshes}
\label{sec:SBPFormulation}
\begin{subequations}
\numberwithin{equation}{section}
In this section, we present a second entropy stable DG method for the SWE based on summation-by-parts (SBP) operators \cite{chen2017entropy, hesthaven2007nodal}. SBP operators are nodal finite-difference matrices which mimic integration by parts. This property is crucial in constructing entropy stable discretizations of non-linear conservation laws. 

The SBP operators used in entropy stable shceme usually include a diagonal mass matrix and differentiation matrices, which are designed to approximate spatial derivatives up to a specified order of accuracy. They are constructed algebraically given a set of quadrature nodes with positive weights and some specific level of accuracy \cite{chen2017entropy, hicken2013summation}. The SBP property is typically sufﬁcient to prove stability for linear PDEs under periodic boundary conditions and appropriate coupling terms \cite{fernandez2014review}. 

Entropy stable numerical schemes can be constructed using either hybridized SBP operators or traditional SBP operators. Hybridized SBP operators provide a numerical connection to polynomial DG methods and allow for flexibility in the choice of quadrature rules. In contrast, traditional SBP operators do not correspond to any polynomial DG discretizations. This is because SBP operators identify nodal values as individual degrees of freedom, and the number of points is larger than the dimension of underlying polynomial space. As a result, traditional SBP operators cannot be derived from standard DG variational formulations. 

SBP operators minic the structure of mass-lumped under-integrated DG discretizations, which can result in higher aliasing errors. However, entropy stable discretizations using traditional SBP operators are more efficient, as they typically involve fewer operations and can skip the entropy projection step discussed in Section \ref{sec:DG_2dSWE}. We will compare the performances of entropy stable methods under both hybridized and traditional SBP operators in Section \ref{sec:numerical_results}.


\subsection{SBP quadrature rules}

We require an SBP-quadrature rule to have the following properties:
\begin{itemize}
\item The surface quadrature points are identically distributed on each face (which enables straightforward coupling between elements.
\item The quadrature weights are positive.
\item The volume quadrature rule is exact for polynomials up to degree $2N-1$.
\item The volume quadrature rule also contains boundary points which form a separate surface quadrature rule. 
\item The surface quadrature rule is exact for degree $2N$ polynomials on each face.
\end{itemize}
In our numerical experiments, we consider two sets of 2D SBP quadrature points. The first uses 1D Gauss-Legendre quadrature on the edges while the second uses 1D Gauss-Lobatto quadrature on edges. They are shown in the Figure \ref{fig:SBP_GLE_nodes} and \ref{fig:SBP_GLO_nodes} respectively. The Gauss-Legendre SBP quadrature rules were introduced in \cite{chen2017entropy}, and the Gauss-Lobatto SBP quadrature rules were computed specifically for this paper.
\begin{figure}[H]
\begin{center}
\includegraphics[width=.23\textwidth]{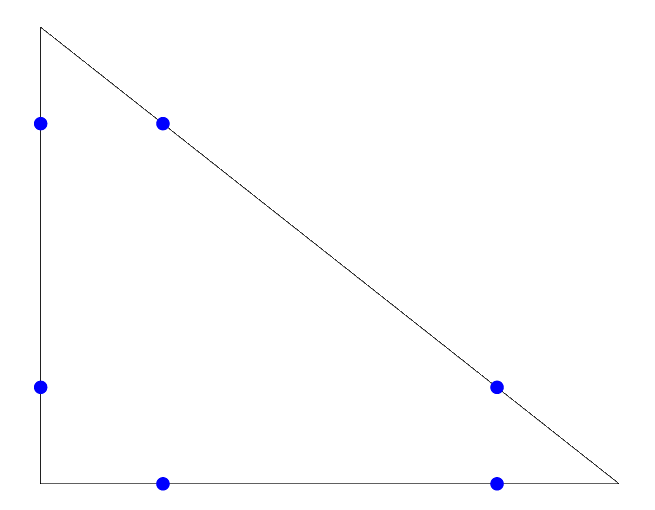}
\includegraphics[width=.23\textwidth]{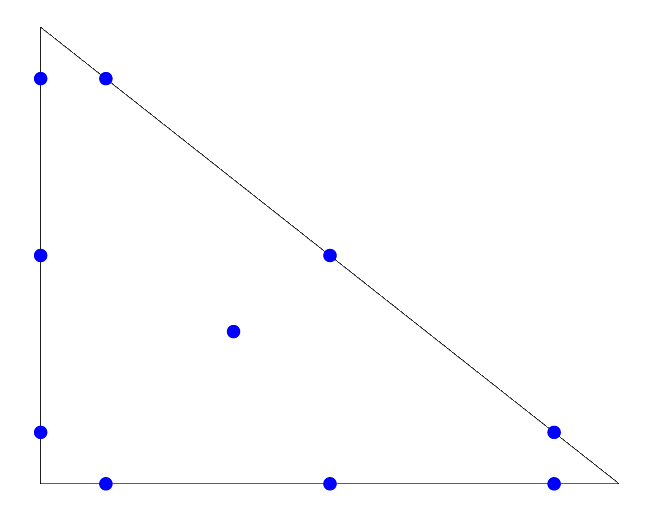}
\includegraphics[width=.23\textwidth]{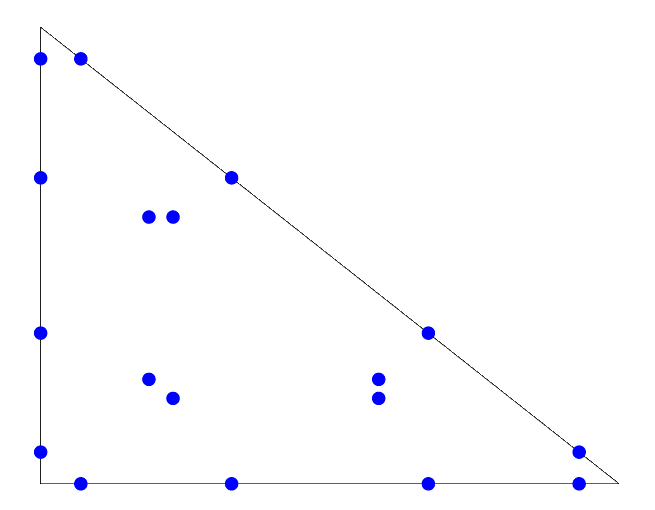}
\includegraphics[width=.23\textwidth]{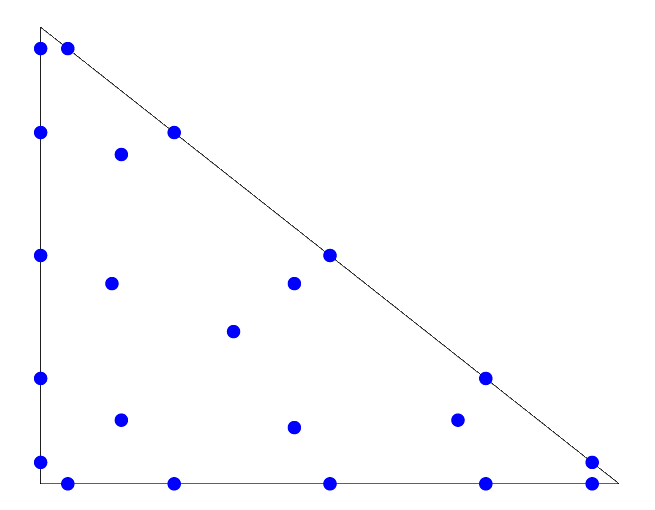}
\caption{Gauss-Legendre quadrature for $N = 1,2,3,4$}
\label{fig:SBP_GLE_nodes}
\end{center}
\end{figure}
\begin{figure}[H]
\begin{center}
\includegraphics[width=.23\textwidth]{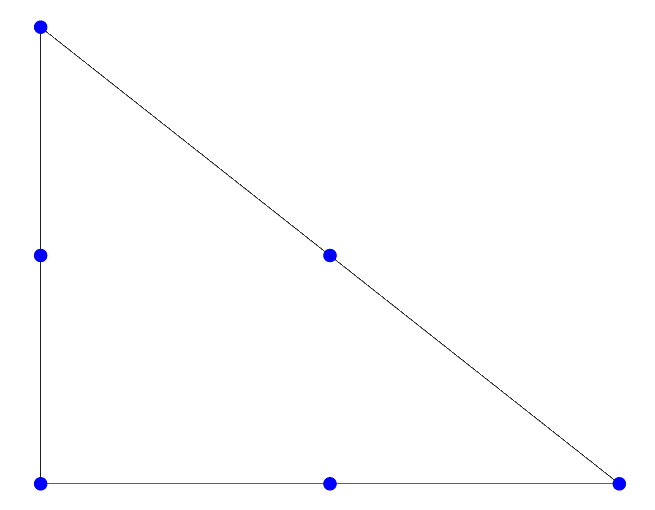}
\includegraphics[width=.23\textwidth]{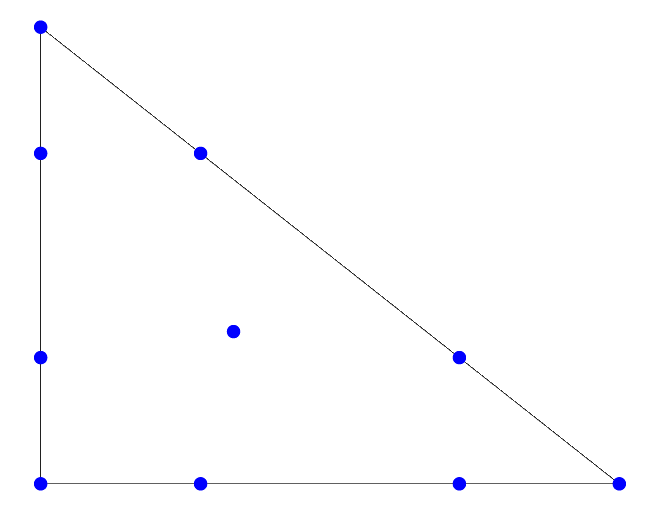}
\includegraphics[width=.23\textwidth]{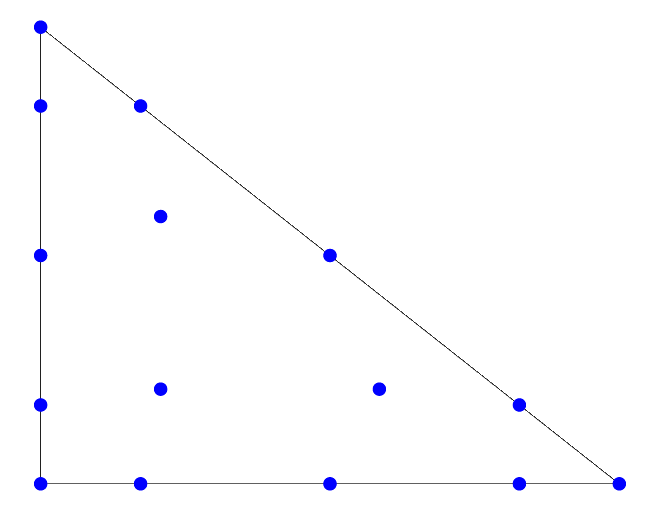}
\includegraphics[width=.23\textwidth]{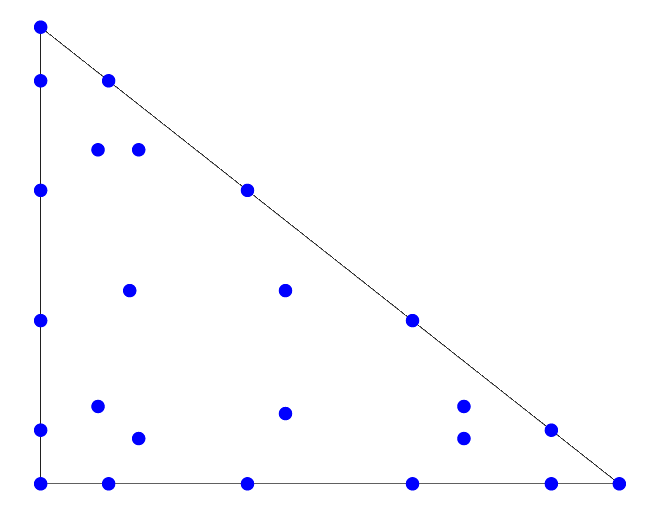}
\caption{Gauss-Lobatto quadrature for $N = 1,2,3,4$}
\label{fig:SBP_GLO_nodes}
\end{center}
\end{figure}

\subsection{Construction of traditional SBP operators}
We now describe a process for constructing traditional SBP operators. We assume we are given a degree $2N-1$ quadrature rule for the polynomial space $P^N(\hat{D})$, with $N_q$ nodes, $\{\bm{x_i}\}_{i=1}^{N_q}$, and positive weights $\{w_i\}_{i=1}^{N_q}$. The nodal values of the function $u(x)$ on the quadrature points is denoted by:
\begin{align}
\bm{u} = [u(\bm{x}_1),  ... , u(\bm{x}_{N_q})]^T.
\end{align}
The SBP mass matrix is defined as the a diagonal matrix with quadrature weights on the diagonal.
\begin{align}
\bm{M}_{\rm{SBP}} = \diag{[w_1, ..., w_{N_q}]}.
\end{align}
We define $\bm{D}^x$ and $\bm{D}^y$ to be the nodal differentiation matrices associated with the $x$ and $y$ derivatives. We also define $\bm{B}^x$ and $\bm{B}^y$ to be diagonal surface matrices, whose entries are surface quadrature weights scaled by the $x$ and $y$ components of the outward normal vector.

We have the following definition \cite{chen2017entropy}.
\begin{Def}{}
Consider the diagonal mass matrix consisting of quadrature weights
\begin{align}
\bm{M}_{\rm{SBP}} = \diag{[w_1, ..., w_{N_q}]}.
\end{align}
A 2D operator $\bm{Q}^i_{\rm{SBP}}$ is said to have SBP property if for $i=x,y$, the following properties holds.
\begin{itemize}
\item Let $\bm{D}^i = \bm{M}_{\rm{SBP}}^{-1}\bm{Q}^i_{\rm{SBP}}$. Then $
(\bm{D}^i\bm{u})_j=\frac{\partial \bm{u}}{\partial x_i}\bigg\rvert_{x = x_j}
$ for any $\bm{u}\in P^N(\hat{D})$.
\item $\bm{Q}^i_{\rm{SBP}}+(\bm{Q}^i_{\rm{SBP}})^T = \bm{B}^i.$
\end{itemize}
\end{Def}
While it is not immediately apparent, one can derive traditional SBP operators from hybridized SBP operators. We introduce the ``selection" matrix $\bm{I}_f$ of size $N_f\times N_q$. $\bm{I}_f$ is a generalized permutation matrix which extracts the surface nodes from the list of all quadrature nodes. For example, suppose that the $i$th node in the list of all quadrature nodes is on the surface of the reference domain. Suppose that this node corrsponds to the $j$th node in the list of surface quadrature nodes. Then, the $(i,j)$ entry of $\bm{I}_f$ is one. To summarize, $\bm{I}_f$ selects out the surface quadrature nodes from the set of SBP quadrature nodes and reorders them for computation.

\begin{Thm}
\label{thm:hybridtoSBP}
Define $\bm{Q}^i_{\rm{SBP}}$ as:
\begin{align}
\bm{Q}^i_{\rm{SBP}} = \begin{bmatrix}
\bm{I}\\\bm{I}_f
\end{bmatrix}^T
\bm{Q}_h^i
\begin{bmatrix}
\bm{I}\\\bm{I}_f
\end{bmatrix},
\end{align}
where $\bm{Q}_h^i$ is defined in Section \ref{subsec:Hybridized_SBP_operators}. Then,  $\bm{Q}^i_{\rm{SBP}}$ a multi-dimensional SBP operator, which is explicitly defined as.
\end{Thm}
\begin{proof}
\label{pf:hybridtoSBP}
In order to show that the matrix $\bm{Q}^i_{\rm{SBP}}$ is consistent with definition 4.1 in \cite{chen2017entropy}, we first need to show that the difference matrix $\bm{M}^{-1}_{\rm{SBP}}\bm{Q}^i_{\rm{SBP}}$ is exact for any polynomial $u(\bm{x}) \in P^N(\hat{D})$. Let $\bm{u}$ denote the modal coefficient of the polynomial and $\bm{u}_q = \bm{V}_q\bm{u}$ denote its nodal value at the quadrature points. We then have
\begin{align}
\bm{Q}^i_{\rm{SBP}} &= \begin{bmatrix}
\bm{I}\\\bm{I}_f
\end{bmatrix}^T
\bm{Q}_h^i
\begin{bmatrix}
\bm{I}\\\bm{I}_f
\end{bmatrix}\\
&=  \begin{bmatrix}
\bm{I}\\\bm{I}_f
\end{bmatrix}^T
\frac{1}{2}
\begin{bmatrix}
\bm{Q}^i - (\bm{Q}^i)^T & \bm{E}^T\bm{B}^i\\
-\bm{B}^i\bm{E} & \bm{B}^i\\
\end{bmatrix}
\begin{bmatrix}
\bm{I}\\\bm{I}_f
\end{bmatrix}.
\end{align}
We use the fact that $\bm{Q}^i + (\bm{Q}^i)^T = \bm{E}^T\bm{B}^i\bm{E}$ from Section \ref{sec:DG_2dSWE} to substitute $(\bm{Q}^i)^T = \bm{E}^T\bm{B}^i\bm{E} - \bm{Q}^i$ into $\bm{Q}^i - (\bm{Q}^i)^T$. We have\begin{align}
\bm{Q}^i - (\bm{Q}^i)^T &= \bm{Q}^i - (\bm{E}^T\bm{B}^i\bm{E} - \bm{Q}^i)\\ 
&= 2\bm{Q}^i - \bm{E}^T\bm{B}^i\bm{E}.
\end{align}
Multiply the $\frac{1}{2}$ into the matrix, we have:
\begin{align}
\bm{Q}^i_{\rm{SBP}} = \begin{bmatrix}
\bm{I}\\\bm{I}_f
\end{bmatrix}^T
\begin{bmatrix}
\bm{Q}^i - \frac{1}{2}\bm{E}^T\bm{B}^i\bm{E} & \frac{1}{2}\bm{E}^T\bm{B}^i\\
-\frac{1}{2}\bm{B}^i\bm{E} & \frac{1}{2}\bm{B}^i\\
\end{bmatrix}
\begin{bmatrix}
\bm{I}\\\bm{I}_f
\end{bmatrix}.
\end{align}
So we can write the difference matrix as:
\begin{align}
\bm{D}^i_{\rm{SBP}} &= \bm{M}^{-1}_{\rm{SBP}}\bm{Q}^i_{\rm{SBP}} \\
&= \bm{M}^{-1}_{\rm{SBP}}\begin{bmatrix}
\bm{I}\\\bm{I}_f
\end{bmatrix}^T
\begin{bmatrix}
\bm{Q}^i - \frac{1}{2}\bm{E}^T\bm{B}^i\bm{E} & \frac{1}{2}\bm{E}^T\bm{B}^i\\
-\frac{1}{2}\bm{B}^i\bm{E} & \frac{1}{2}\bm{B}^i\\
\end{bmatrix}
\begin{bmatrix}
\bm{I}\\\bm{I}_f
\end{bmatrix}.
\end{align}
If we apply $\bm{D}^i_{\rm{SBP}}$ at $\bm{u}_q$, we have
\begin{align}
\bm{D}^i_{\rm{SBP}}\bm{u}_q = \bm{M}^{-1}_{\rm{SBP}}\begin{bmatrix}
\bm{Q}^i - \frac{1}{2}\bm{E}^T\bm{B}^i\bm{E} -\frac{1}{2}\bm{I}_f^T\bm{B}^i\bm{E},& \frac{1}{2}\bm{E}^T\bm{B}^i + \frac{1}{2}\bm{I}_f^T\bm{B}^i\\
\end{bmatrix}
\begin{bmatrix}
\bm{u}_q \\ \bm{u}_f
\end{bmatrix},
\end{align}
since $\begin{bmatrix}\bm{I}\\\bm{I}_f\end{bmatrix}$ maps the vector $\bm{u}_q$ to the vector of values at both volume and surface quadrature points. Notice that we have $\bm{E} = \bm{V}_f\bm{P}_q$ and $\bm{P}_q\bm{V}_q = \bm{I}$ from \ref{sec:DG_2dSWE}. Then, $\bm{E}\bm{u}_q = \bm{V}_f\bm{P}_q\bm{V}_q\bm{u} = \bm{V}_f\bm{u} = \bm{u} _f$

\begin{align}
\nonumber\bm{D}^i_{\rm{SBP}}\bm{u}_q &= \bm{M}^{-1}_{\rm{SBP}}\LRp{\bm{Q}^i\bm{u}_q - \frac{1}{2}\bm{E}^T\bm{B}^i\bm{E}\bm{u}_q -\frac{1}{2}\bm{I}_f^T\bm{B}^i\bm{E}\bm{u}_q+\frac{1}{2}\bm{E}^T\bm{B}^i\bm{u}_f + \frac{1}{2}\bm{I}_f^T\bm{B}^i\bm{u}_f}\\
 &=\bm{M}^{-1}_{\rm{SBP}}\LRp{\bm{Q}^i\bm{u}_q - \frac{1}{2}\bm{E}^T\bm{B}^i\bm{u}_f -\frac{1}{2}\bm{I}_f^T\bm{B}^i\bm{u}_f+\frac{1}{2}\bm{E}^T\bm{B}^i\bm{u}_f + \frac{1}{2}\bm{I}_f^T\bm{B}^i\bm{u}_f}.
\end{align}
We simplify the above formulation and substitute in $\bm{Q}^i = \bm{P}_q^T\bm{\hat{Q}}^i\bm{P}_q$, $\bm{P}_q = \bm{M}^{-1}\bm{V}_q^T\bm{W}$ and $\bm{M}^{-1}_{\rm{SBP}} = \bm{W}$. Since $\bm{M}$ is symmetric, $\bm{M}^{-1} = (\bm{M}^{-1})^T$. Then using $\bm{P}_q\bm{V}_q = \bm{I}$,
\begin{align}
\bm{D}^i_{\rm{SBP}}\bm{u}_q &= \bm{M}^{-1}_{\rm{SBP}}\bm{Q}^i\bm{V}_q\bm{u}\\
&= \bm{W}^{-1}\bm{P}_q^T\bm{\hat{Q}}^i\bm{P}_q\bm{V}_q\bm{u}\\
&= \bm{W}^{-1}\bm{W}\bm{V}_q(\bm{M}^{-1})^T\bm{M}\bm{D}^i(\bm{P}_q\bm{V}_q)\bm{u}\\
&= \bm{V}_q\bm{D}^i\bm{u}.
\end{align}
Recall that the differentiation matrix $\bm{D}^i$ is exact for polynomials $\bm{u}\in P^N(\hat{D})$. Since $\bm{V}_q$ is a degree $N$ interpolation matrix of the derivative at quadrature points, we conclude that $\bm{D}^i_{\rm{SBP}}$ exactly differentiate $\bm{u}_q$, where $\bm{u}_q$ is the vector of values of a polynomial $\bm{u}\in P^N(\hat{D})$ at quadrature points.

We now show the summation-by-parts property. First we use the property that
\begin{align}
\bm{Q}^i_h + (\bm{Q}^i_h)^T &= \begin{bmatrix}
\bm{0} & \\
  &\bm{B}^i \\
\end{bmatrix},
\end{align}
to rewrite $\bm{Q}^i_{\rm{SBP}}$ as
\begin{align}
\bm{Q}^i_{\rm{SBP}} &= \begin{bmatrix}
\bm{I}\\\bm{I}_f
\end{bmatrix}^T
\LRp{\begin{bmatrix}
\bm{0} & \\
  &\bm{B}^i \\
\end{bmatrix} - (\bm{Q}_h^i)^T}
\begin{bmatrix}
\bm{I}\\\bm{I}_f
\end{bmatrix},\\
&= \begin{bmatrix}
\bm{I}\\\bm{I}_f
\end{bmatrix}^T
\begin{bmatrix}
\bm{0} & \\
  &\bm{B}^i \\
\end{bmatrix}
\begin{bmatrix}
\bm{I}\\\bm{I}_f
\end{bmatrix} - \begin{bmatrix}
\bm{I}\\\bm{I}_f
\end{bmatrix}^T
(\bm{Q}_h^i)^T
\begin{bmatrix}
\bm{I}\\\bm{I}_f
\end{bmatrix},\\
&= \bm{I}_f^T\bm{B}^i\bm{I}_f - (\bm{Q}^i_{\rm{SBP}})^T.
\end{align}
We conclude that $\bm{Q}^i_{\rm{SBP}} + (\bm{Q}^i_{\rm{SBP}})^T =\bm{I}_f^T\bm{B}^i\bm{I}_f$.
Since $\bm{I}_f$ is a matrix which selects the surface quadrature points, $\bm{I}_f^T\bm{B}^i\bm{I}_f$ is still diagonal with entries of $\bm{B}^i$ permuted. The permutation order depends on the order of the quadrature points listed in the implementation, but does not affect the summation-by-parts property.
\end{proof}
Given this equivalence, we can now construct an entropy conservative DG-SBP form:
\begin{align}
\bm{M}\td{\bm{u}}{t} + \sum_{i=x,y} \LRp{2\bm{Q}^i_{\rm{SBP}}\circ \bm{F}^i}\bm{1} + \bm{I}_f^T\bm{B}^i\LRp{\bm{f}^i_S-\bm{f}^i} = \bm{S}.
\end{align}
The stability analysis of the SBP formulation follows from the results from \cite{chen2017entropy}. The extension to curved elements can be found in \cite{crean2018entropy}. 
\end{subequations}
\section{Numerical results}
\label{sec:numerical_results}
\begin{subequations}
\numberwithin{equation}{section}
In this section, we present some two dimensional numerical experiments and results to demonstrate the accuracy and stability of the entropy stable DG scheme. All experiments are run using an entropy stable scheme, which we construct by adding local Lax-Friedrichs penalization \cite{wintermeyer2017entropy} to a baseline entropy conservative DG formulation.

The first experiment is a ``lake-at-rest" condition to test the well-balancedness of our scheme. The second experiment is a translating vortex. This problem has an explicit analytic solution, which we use to investigate the convergence rate of our algorithm. The third experiment is a dam break simulation from \cite{wintermeyer2017entropy}. The last experiment is a converging channel simulation \cite{wirasaet2015artificial}. 

All numerical experiments utilize the fourth order five-stage low-storage Runge-–Kutta method \cite{carpenter1994fourth}. Following the derivation of stable timestep restrictions in \cite{chan2016gpu}, we define the timestep $\Delta t$ to be
\begin{align}
\Delta t= CFL \times \frac{h}{C_N},\hspace{1cm} C_N = \frac{(N+1)(N+2)}{2},
\end{align}
where $C_N$ is the degree dependent constant in the trace inequality space \cite{warburton2003constants}, and CFL is a user-defined constant. We use CFL = 0.125 for all experiments. 

We test ``lake-at-rest" and translating vortex problem on both affine and curved meshes. For the curved mesh, we construct the warping of the regular triangular mesh in the following way:
\begin{flalign*}
&x = x + C_{curve}L_x\cos\LRp{\pi\LRp{x-x_0}/L_x}\cos\LRp{1.5\pi\LRp{y-y_0}/L_y},\\
&y = y + C_{curve}L_y\sin\LRp{2\pi\LRp{x-x0}/L_x}\cos\LRp{pi\LRp{y-y0}/L_y},
\end{flalign*}
where $(x_0,y_0) = (0,0)$ is location of the center of the simulation, $(0,0)$ in our case. $L_x$ and $L_y$ are the length of the domain in $x$ and $y$ direction respectively. $C_{curve}$ is a user specified curving coefficient, which we use 0.1 in all the experiments in "lake-at-rest" and translating vortex problem. For the dam break and converging channel experiment, we construct a curved mesh which is boundary-fitted to the curved dam (a curved interior boundary) or the curved channel. More details are provided in the following sections. 

\subsection{Lake at rest}
We first consider the ``lake-at-rest" condition \cite{gassner2016well, noelle2007high, leveque1998balancing} on $[-1,1]\times[-1,1]$ and we set the boundary to be periodic. We set:
\begin{align}
H &= h + b  = 2,\\
b(x,y) &= 0.1\sin(2\pi x)\cos(2\pi x) + 0.5.
\end{align}

For a degree N approximation, we first interpolate $b$ using a continuous $C^0$ degree N polynomial, then we set $h=2-b$. In the Table \ref{Tab:WB_affine} and Table \ref{Tab:WB_curved}, we observe the error for this setting is of the magnitude of round-off errors at each quadrature point for both affine and curved triangular meshes. All of the ``lake-at-rest" experiments are run to a final time of $T = 1/2$. 
\begin{table}[H]

\begin{center}
\begin{tabular}{|c|c|c|c|c|c|}
\hline
   $N$  &  $K$ = 256  &  $K$ = 1024  &  $K$ = 4096  &  $K$ = 16384 \\
\hline
   1  &  1.03E-13 &  1.54E-13  &  3.17E-13  &  7.27E-13  \\
\hline
   2  &  9.62E-13 &  4.61E-12  &  1.90E-11  &  7.88E-11  \\
\hline
   3  &  2.96E-12 &  7.71E-12  &  4.43E-11  &  1.58E-10  \\
\hline
   4  &  1.62E-11 &  6.41E-11  &  2.53E-10  &  1.03E-09  \\  
\hline
\end{tabular}
\end{center}
\caption{$L^2$ error for the lake at rest problem on affine triangular meshes. $N$ denotes the polynomial degree and $K$ denotes the number of elements.}
\label{Tab:WB_affine}
\end{table}

\begin{table}[H]
\begin{center}
\begin{tabular}{|c|c|c|c|c|c|}
\hline
   $N$  &  $K$ = 256  &  $K$ = 1024  &  $K$ = 4096  & $K$ = 16384 \\
\hline
   1  &  1.06E-13 & 1.64E-13 & 2.88E-13 & 5.51E-13\\
\hline
   2  &  1.07E-12 & 4.06E-12 & 1.67E-11 & 6.78E-11 \\
\hline
   3  &  2.99E-12 & 1.16E-11 & 4.55E-11 & 1.84E-10 \\
\hline
   4  &  1.17E-11 & 4.66E-11 & 1.90E-10 & 7.57E-10\\
\hline
\end{tabular}
\end{center}
\caption{$L^2$ error for the lake at rest problem on curved triangular mesh. $N$ denotes the polynomial degree and $K$ denotes the number of elements.}
\label{Tab:WB_curved}
\end{table}
We evaluate the $L^2$ error using a more accurate triangular quadrature rule exact for degree $2N+2$ polynomials. We notice that the error for this problem increases as we use higher order polynomial bases or finer meshes, and we attribute this to numerical round off.

\subsection{Translating vortex}
We now consider the vortex translation test. We set the domain to be $[-10,10]\times[-5,5]$ and the exact solution for the vortex at any time $t$ is given by \cite{kang2020imex, ricchiuto2009stabilized}:
\begin{align}
&h = h_{\infty}-\frac{\beta^2}{32\pi^2}e^{-2(r^2-1)},\\
&u = u_{\infty}-\frac{\beta}{2\pi}e^{-2(r^2-1)}y_t,\\
&v = v_{\infty}+\frac{\beta}{2\pi}e^{-2(r^2-1)}x_t,\\
b = 0,
\end{align}
where
\begin{align}
x_t &= x - x_c - u_{\infty}t,\\
y_t &= y - y_c - v_{\infty}t,\\
r^2 &= x^2 + y^2.
\end{align}
In this example, we set 
\begin{align}
h_{\infty} = 1,\quad \beta = 5, \quad  g = 2 \quad\rm{and}\quad (u_{\infty}, v_{\infty}) = (1,0).
\end{align}
Initially the vortex is located at $(x_c, y_c) = (0, 0)$. In this setup, the vortex propagates to the right along the x-axis.  The domain and problem setup are chosen such that periodic boundary conditions can be used without affecting accuracy.

We use both affine and curved meshes for this experiment. We also compare the $L^2$ error the SBP formulation using operators based on both Gauss-Legendre and Gauss-Lobatto quadrature nodes as we introduced in section \ref{sec:SBPFormulation}. We calculate the $L^2$ error using the same method as in the ``lake-at-rest" problem. 
\begin{figure}
\begin{center}
\includegraphics[width=.75\textwidth]{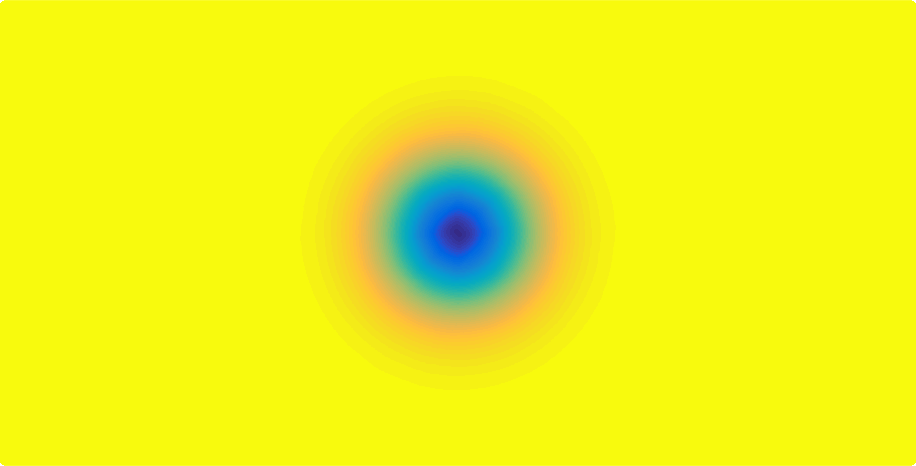}
\caption{A translating vortex in 2D}
\label{fig:vortex}
\end{center}
\end{figure}
\begin{figure}
\begin{center}
\includegraphics[width=.45\textwidth]{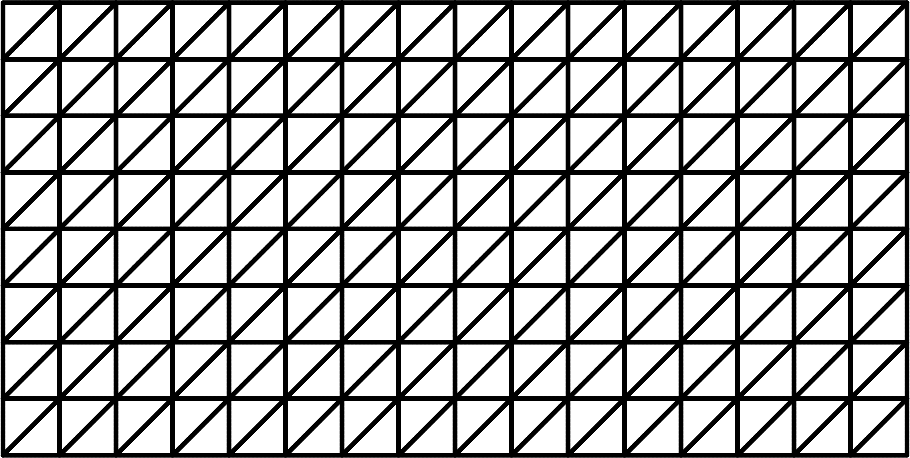}
\includegraphics[width=.45\textwidth]{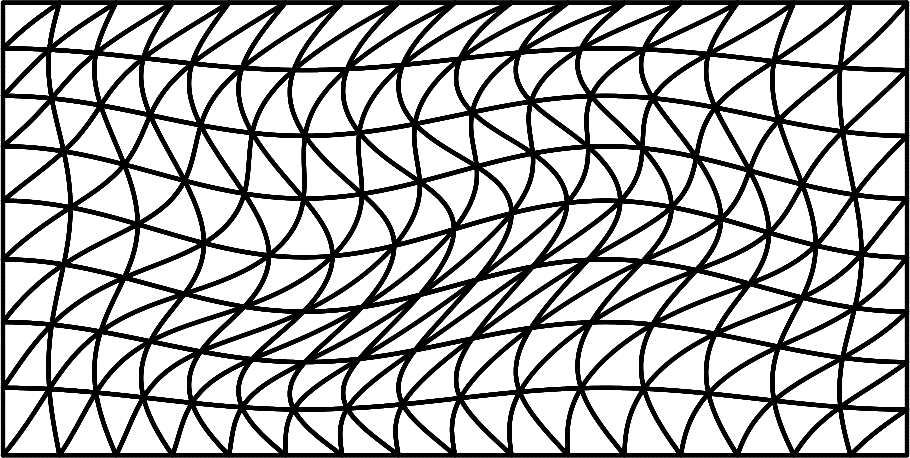}
\caption{Meshes in 2D for translating vortex problem}
\label{fig:mesh vortex}
\end{center}
\end{figure}
 The convergence results are presented in the Figure \ref{fig:vortex_error1} and Figure \ref{fig:vortex_error2}:

\begin{figure}
\begin{center}
\includegraphics[width=0.8\textwidth]{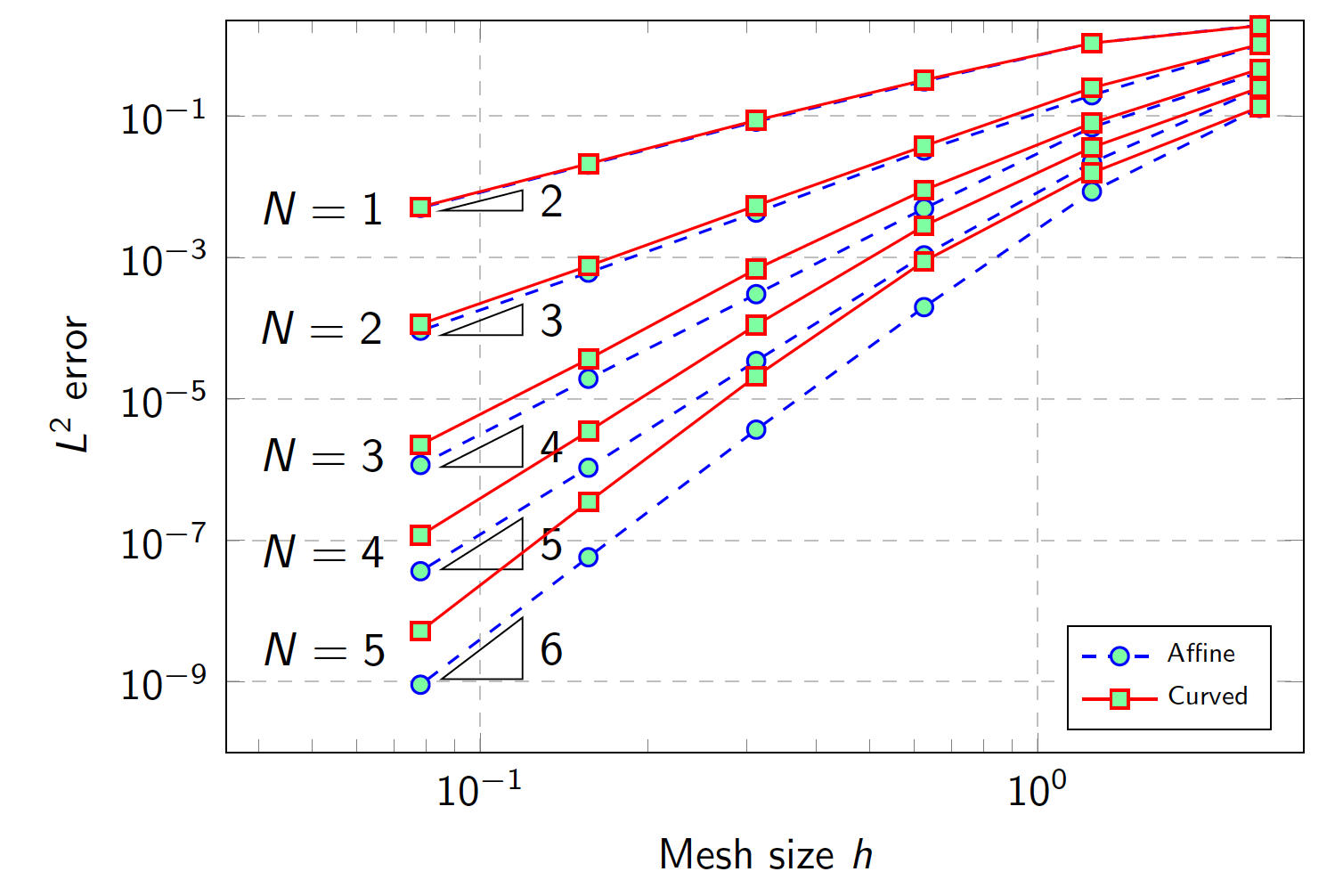}
\caption{$L_2$ error for the translating vortex after 0.5 second on affine and curved meshes}
\label{fig:vortex_error1}
\end{center}
\end{figure}

\begin{figure}
\begin{center}
\includegraphics[width=0.8\textwidth]{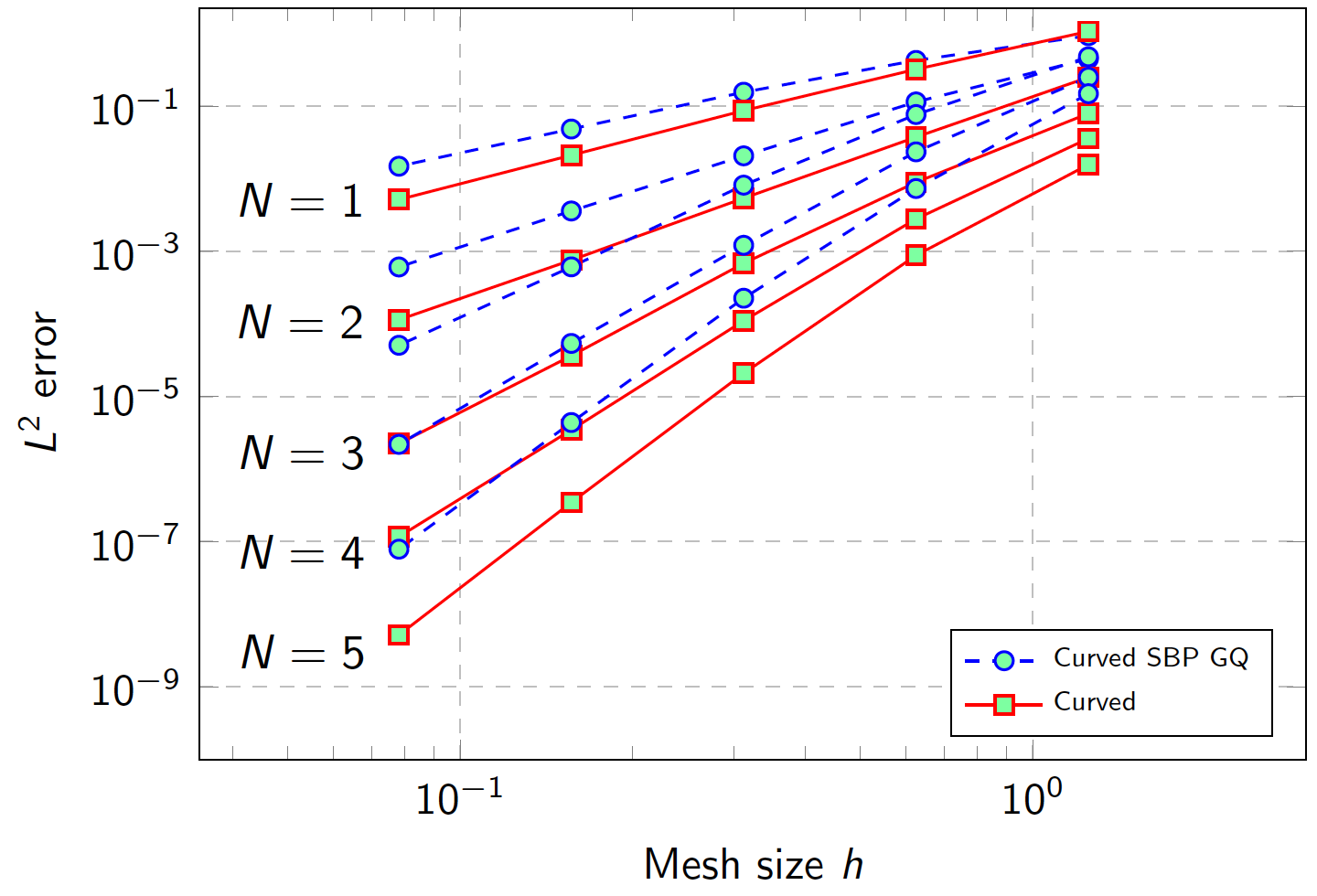}
\caption{$L_2$ errors for the translating vortex after 0.5 seconds using hybridized DG and Gauss-Legendre SBP operator schemes on curved meshes.}
\label{fig:vortex_error2}
\end{center}
\end{figure}

\begin{figure}
\begin{center}
\includegraphics[width=0.8\textwidth]{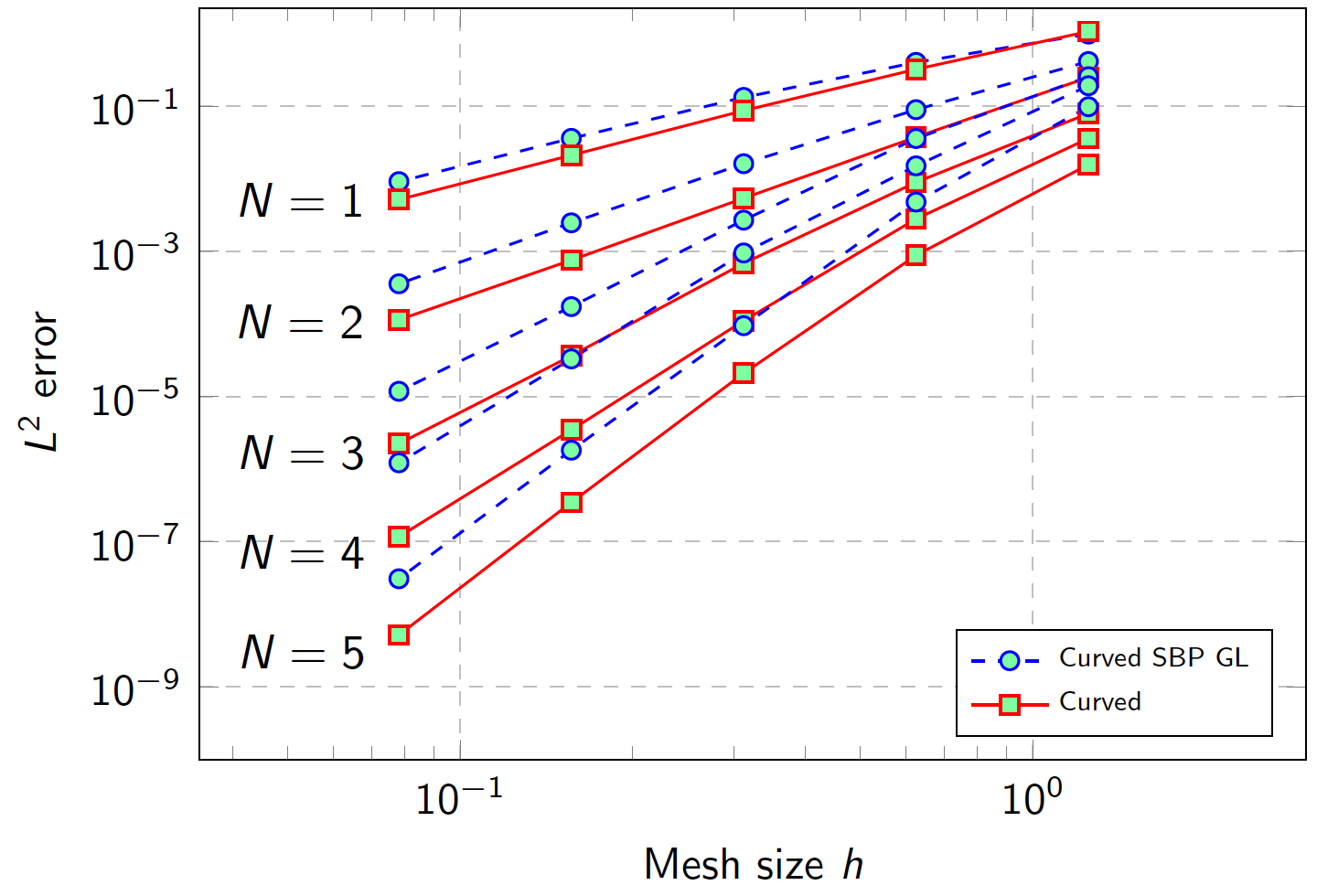}
\caption{$L_2$ errors for the translating vortex after 0.5 seconds using hybridized DG and Gauss-Lobatto SBP operator schemes on curved meshes.}
\label{fig:vortex_error3}
\end{center}
\end{figure}

Recall that the SBP-DG discretization does not correspond to a polynomial approximation space.  Thus, to calculate the $L^2$ error for the SBP-DG discretization, we first project the final numerical solution to polynomials of degree $N$. We then use this projection to evaluate the $L^2$ error using a quadrature rule which is exact for at least degree $2N+2$ polynomials. The error for the DG method with hybridized SBP operators is computed using the same quadrature. 

In Figure \ref{fig:vortex_error1}, we analyze the accuracy of the hybridized SBP scheme on both curved and affine meshes. We observe the same rate of convergence for both curved and affine meshes, but note that some accuracy is lost on the curved meshes. In Figure \ref{fig:vortex_error2} and \ref{fig:vortex_error3}, we show the comparisons between two SBP methods with the hybridized SBP method. Figure \ref{fig:vortex_error2} presents the Gauss-Legendre SBP nodes and Figure \ref{fig:vortex_error3} presents the Gauss-Lobatto nodes compared against the hybridized SBP method. We observe that the use of traditional SBP operators, while more efficient, lose one order of accuracy compared to the use of hybridized SBP operators. The Gauss-Lobatto SBP operator seems to be slightly more accurate than Gauss-Lobatto SBP operator, but both converge at about the same rate.

\subsection{Dam break}
This experiment is taken from \cite{wintermeyer2018entropy, lukacova2009entropy}. We utilize the same physical setting but use curved triangular meshes instead of curved quadrilateral meshes. The domain $[-10,10]^2$ is discretized using a $20\times20$ grid of quadrilaterals, which are then split into triangular elements. We set $N=3$ for the polynomial degree. 

The dam is modeled by imposing reflective boundary conditions along the curve defined by the following function:
\begin{align}
x = \frac{1}{25}y^2,
\end{align}
with a break between $y=-0.5$ and $y=0.5$ to allow water to flow through, as shown in Figure \ref{fig:DamMesh}. The dam is marked in red and fitted by the curved mesh. 

We start with a constant water height on both sides of the dam with zero initial velocity. The bathymetry is set to $b = 0$ on both sides. We set the initial water height to be $h = 10$ on the left side of the dam and $h=5$ on the right side as shown in Figure \ref{fig:DamMesh}.
\begin{figure}[H]
\begin{center}
\includegraphics[width=0.45\textwidth]{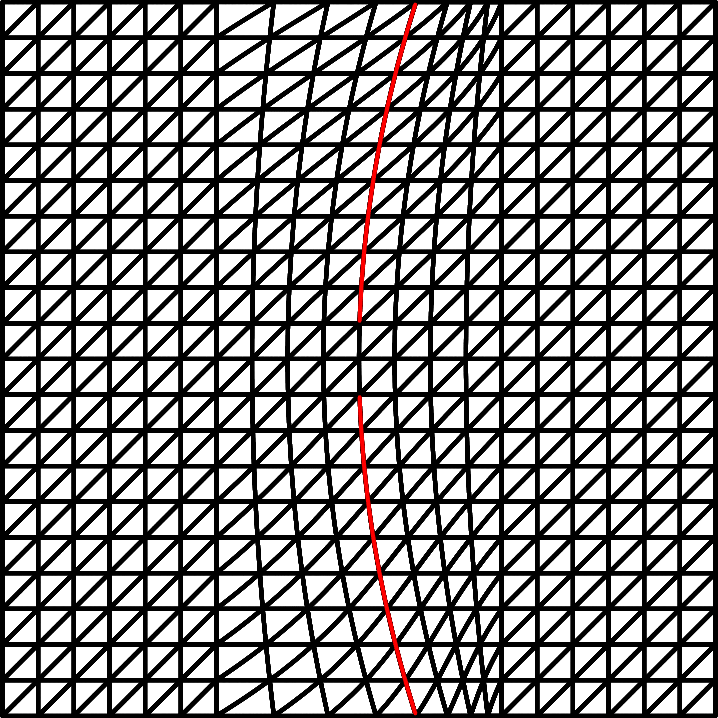}
\includegraphics[width=0.5\textwidth]{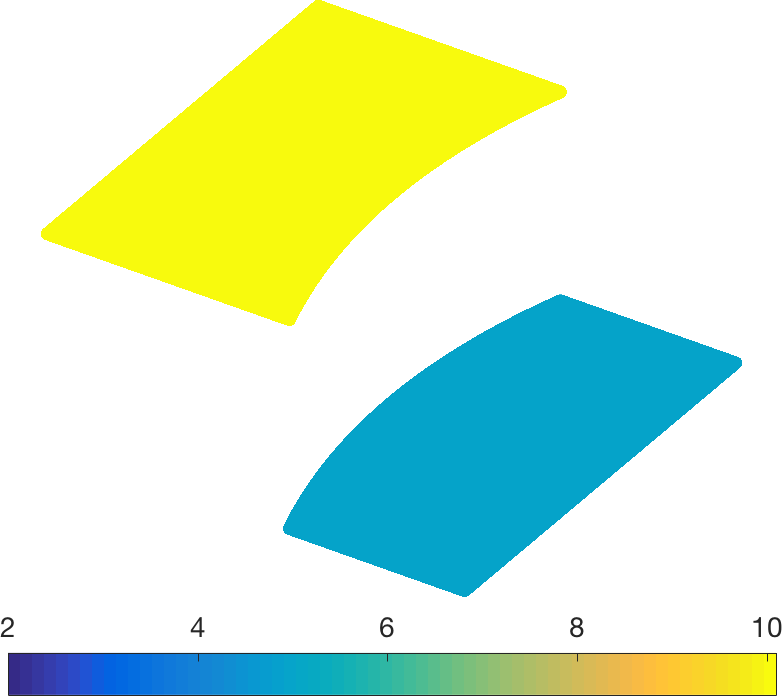}\\
\caption{left: curved mesh for 2D dam break problem with the red curve denoting the broken dam. right: initial condition}
\label{fig:DamMesh}
\end{center}
\end{figure}
We plot for the solution at several different times in Figures \ref{fig:Dam1}, \ref{fig:Dam2} and \ref{fig:Dam3}.We observe that the water falling from the left side of the dam to the right produce a wave front in the lower half. This wave is discontinuous, so we observe some mesh dependent solution oscillations, but they are on a smaller scale and do not cause the solution to blow up. We have also tried simulations with polynomial degrees $N= 5$ and $N=7$. Both of these simulations also remain stable throughout the duration (1.5s) of the run. The numerical solution of the parabolic dam break problem demonstrates that the entropy stable numerical remains robust on the presence of shock discontinuities.
\begin{figure}[H]
\begin{center}
\includegraphics[width=0.45\textwidth]{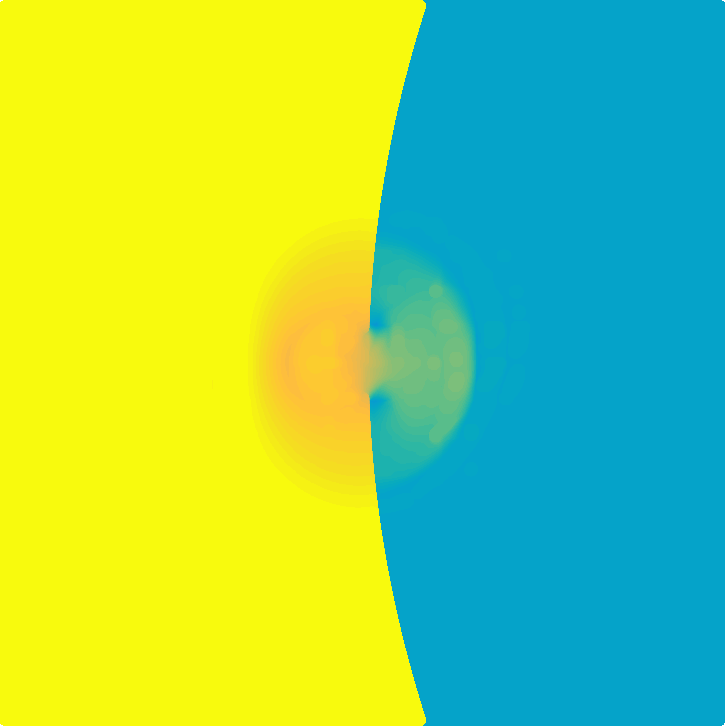}
\includegraphics[width=0.5\textwidth]{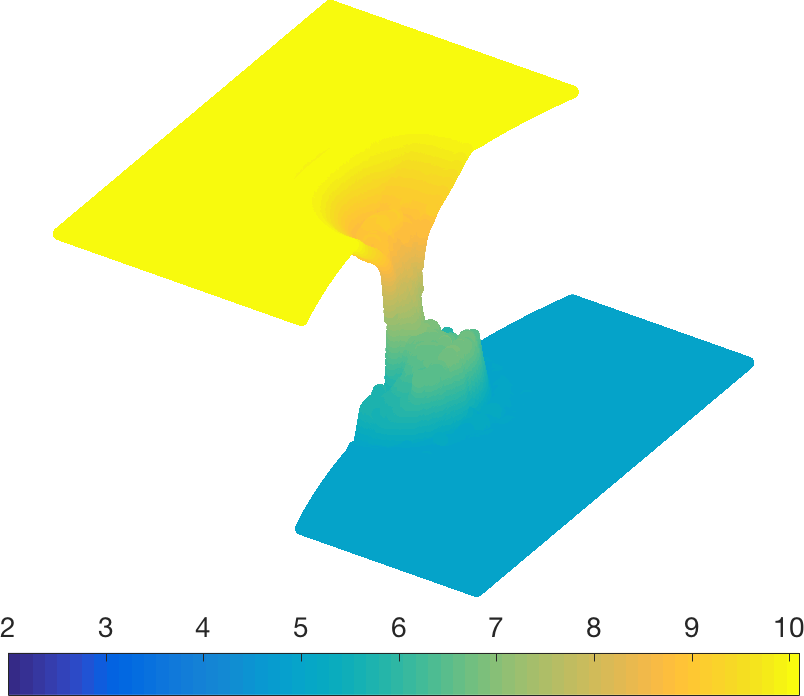}\\
\caption{Top  and side view of the dam breaking problem at T=0.5s}
\label{fig:Dam1}
\end{center}
\end{figure}

\begin{figure}[H]
\begin{center}
\includegraphics[width=0.45\textwidth]{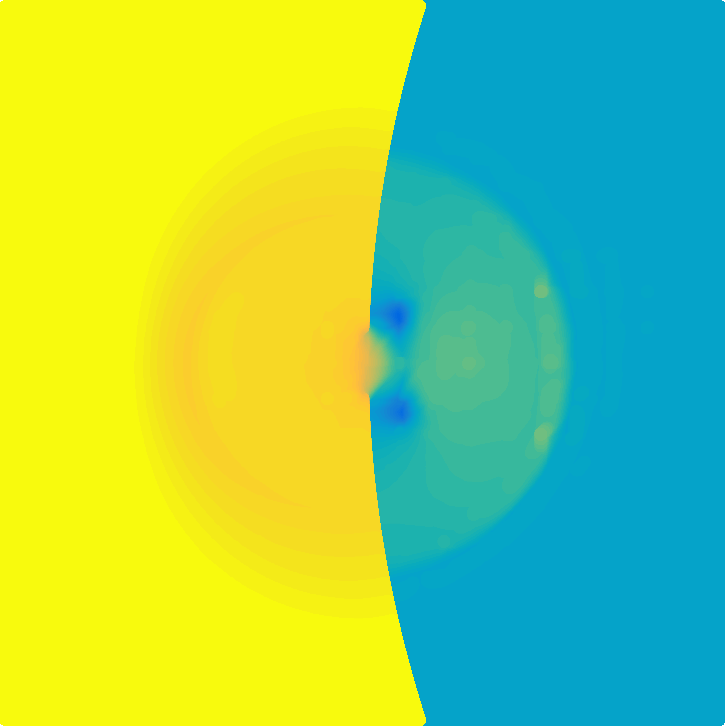}
\includegraphics[width=0.5\textwidth]{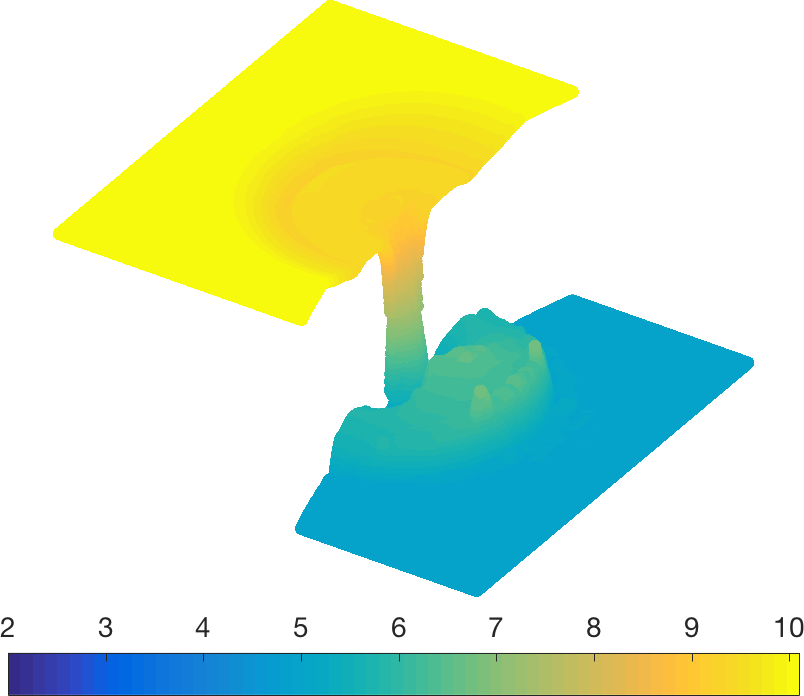}\\
\caption{Top  and side view of the dam breaking problem at T=1s}
\label{fig:Dam2}
\end{center}
\end{figure}

\begin{figure}[H]
\begin{center}
\includegraphics[width=0.45\textwidth]{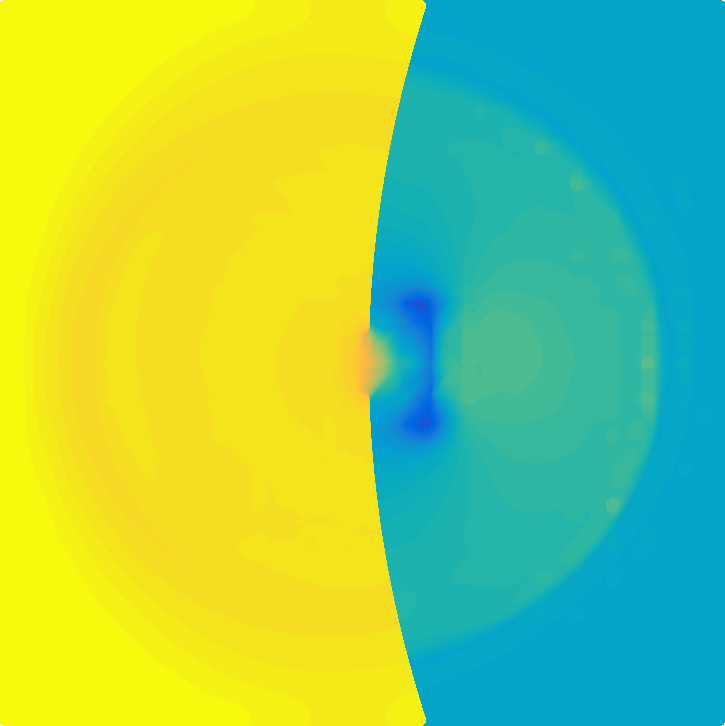}
\includegraphics[width=0.5\textwidth]{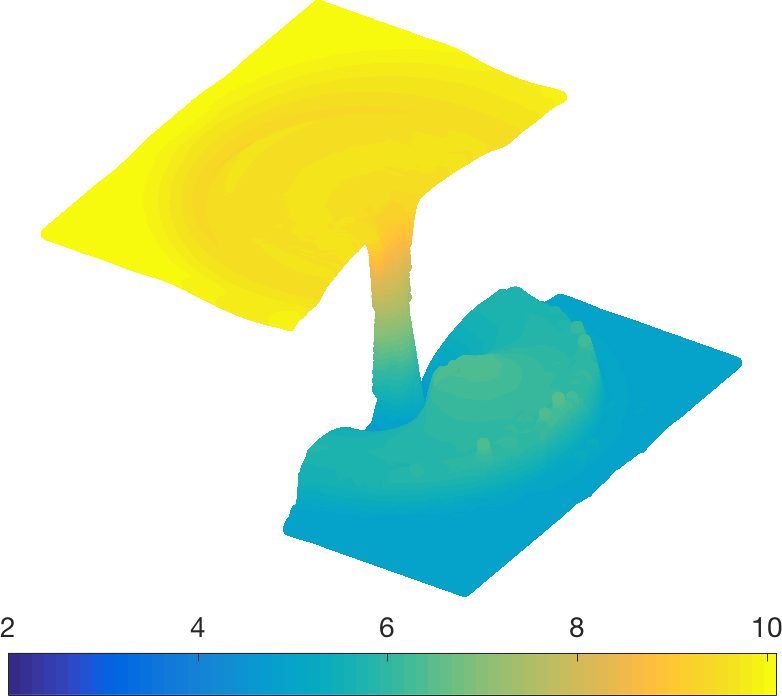}\\
\caption{Top  and side view of the dam breaking problem at T=1.5s}
\label{fig:Dam3}
\end{center}
\end{figure}
\par\bigskip 

\end{subequations}
\section{GPU optimization}
\label{sec:GPU_opt}
\begin{subequations}
\numberwithin{equation}{section}
The use of Graphic Processing Units (GPUs) in scientific computing has become well-established over the last 20 years. For example, the use of GPUs for accelerating the solution of PDE using finite difference methods is common in seismic applications \cite{michea2010accelerating}. GPUs have been widely used to accelerate explicit time-stepping schemes, which typically have high arithmetic intensity \cite{chan2016gpu}. DG algorithms with explicit time integration are well suited to the parallel GPU architectures since most of the computational work is element local, and elements are only loosely coupled through shared interfaces. For each element, the computational intensity is also very high. In this section, we present the details of our GPU implementation and describe some computational optimization applied.

\subsection{Infrastructure}
We start by introducing our general coding infrastructure. The GPU implementation is written with C++ and OCCA code. OCCA is a unified approach to multi-threading languages, which compiles code written in the OCCA kernel language (OKL) at runtime for either CPU (Serial, OpenMP) or GPU (OpenCL, CUDA) architectures. We run all experiments on Google Cloud using a NVIDIA Tesla V100 GPU. We use double precision for all of our tests and the Tesla V100 is a device designed for scientific computations and better suited for double precision calculations than other general purpose GPUs.

We construct the mesh, quadrature nodes, and matrices on a CPU before moving on to the main time-stepping. We divide the iteration into four OCCA kernels: projection, volume, surface and update. The projection kernel performs the entropy projection. The volume kernel calculates volume contributions while the surface kernel accumulates contributions from surface fluxes for each element. Finally, the update kernel applies the inverses of mass matrices (which are precomputed and stored on the GPU) on each element and performs time-integration.

In our entropy stable scheme, we apply flux differencing in the volume kernel, which computes the Hadamard product of differentiation and flux matrices over each element. The entries of the flux matrix are constructed on-the-fly by evaluating entropy conservative flux functions between each pair of nodal values of the solution.  This operation involves significantly more computations than standard matrix vector multiplication, and is expected to be well suited to GPUs architectures. Despite this, the volume kernel is the most expensive kernel in our implementation.

\subsection{Optimization}
Following \cite{wintermeyer2018entropy}, we have applied some optimizations to our GPU implementation:
\begin{itemize}
    \item Step 1: Utilizing Shared Memory.\\
    The first step in improving performance is to reduce the number of reads from  GPU global memory. We load all solution values on an element into the shared memory, which avoids repeated access to slow global memory.
    \item Step 2: Declaring variables constant and pointers restricted.\\
    When a variable does not mutate during its lifespan, we declare the variable constant. We also use the restricted keyword on pointers to memory that is not referenced by other pointers. This allows the compiler to optimize performance. We also pre-compute all constants, (such as $\frac{1}{2}g$), before the start of the kernel to avoid repeated computation in each iteration. 
    \item Step 3: Processing multiple elements per block.\\
    GPUs perform operations in synchronized groups of 32 threads, referred to as ``warps''. However, the natural number of threads used in each block is not typically a multiple of 32, which leads to idle threads and reduced performance on GPUs. In order to better utilize GPU resources and reduce the number of idle threads, we process multiple elements at the same time.  In practice we manually optimize the block size to minimize run time of each kernel.  
    
    This step, however, has its limitations. As the degree of the polynomial basis increases, the amount of shared memory used by each element also increases, but we cannot exceed the total GPU shared memory when combining multiple elements in to the same block. 
    
    \item Step 4: Finally, we utilize an optimized implementation of the volume kernel based on matrix structure.\\
    The first three steps are common and can be applied to general GPU programming. Step 4 is more specific to the structure of the entropy stable DG formulation using hybridized SBP operators. 
    
    We can reduce the computational cost of computing $(\bm{Q}^i_h \circ \bm{F})\bm{1}$, by avoiding computing the Hadamard product for any zero sub-blocks of the matrix $\bm{Q}^i_h$. Using the fact that $\bm{Q}^i_h = \bm{B}^i_h-(\bm{Q}^i_h)^T$, we can rewrite the DG formulation without projection as
    \begin{align}
    \bm{M}\td{\bm{u}}{t} + \sum_{i=x,y} \LRp{(\bm{Q}^i_h - (\bm{Q}^i_h)^T + \bm{B}^i_h)\circ \bm{F}^i}\bm{1} + \bm{B}^i_h\LRp{\bm{f}^i_S(\bm{u}^+,\bm{u})-\bm{f}^i(\bm{u})} = \bm{S}.
    \end{align}
    Using the identity $(\bm{B}^i_h\circ \bm{F}^i)\bm{1}  = \bm{B}^i_h\bm{f}^i(\bm{u})$, which follows from the consistency of the flux, we have
    \begin{align}
    \bm{M}\td{\bm{u}}{t} + \sum_{i=x,y} \LRp{(\bm{Q}^i - (\bm{Q}^i)^T)\circ \bm{F}_i}\bm{1} + \bm{B}^i\LRp{\bm{f}^i_S(\bm{u}^+,\bm{u})} = \bm{S}.
    \end{align}
    We define $\bm{Q}^i_{h,skew} = \bm{Q}^i_h - (\bm{Q}^i_h)^T$. Then, assembling all the pieces and combining the projection step, we can rewrite our DG scheme out as: 
    \begin{align}
    \bm{M}\td{\bm{u}}{t} + \sum_{i=x,y} \begin{bmatrix}
    \bm{V}_q\\\bm{V}_f
    \end{bmatrix}^T\LRp{2\bm{Q}_{h,skew}^i\circ \bm{F}^i}\bm{1} + \bm{V}_f^T\bm{B}^i_h\bm{f}^i_S = \bm{S}.
    \end{align}
    
    We observe that the matrix $\bm{Q}^i_{h,skew}$ has block structure where the lower right block is all zeros as shown in Figure \ref{fig:OCCA4}. We split the original OCCA kernel that computes the Hadamard product of $\bm{Q}^i_{h,skew}$ and $\bm{F_i}$ into two kernels. The first computes the Hadamard product of the first $N_q$ columns and accumulates the row sum into a vector, which corresponds to the part of the matrix in red box in Figure \ref{fig:OCCA4}. The second kernel repeats the same process the upper right part of the matrix in the blue box as in Figure \ref{fig:OCCA4}, then updates the accumulating vector counter for the first $N_q$ entries. This eliminates the work of computing the Hadamard product of zeros entries in the matrix $\bm{Q}^i_{h,skew}$, which consists of approximately one quarter of the computational work.
    
    \begin{center}
    \end{center}
    \vspace{-.5em}
    \begin{figure}
    \centering
    \includegraphics[width=0.75\textwidth]{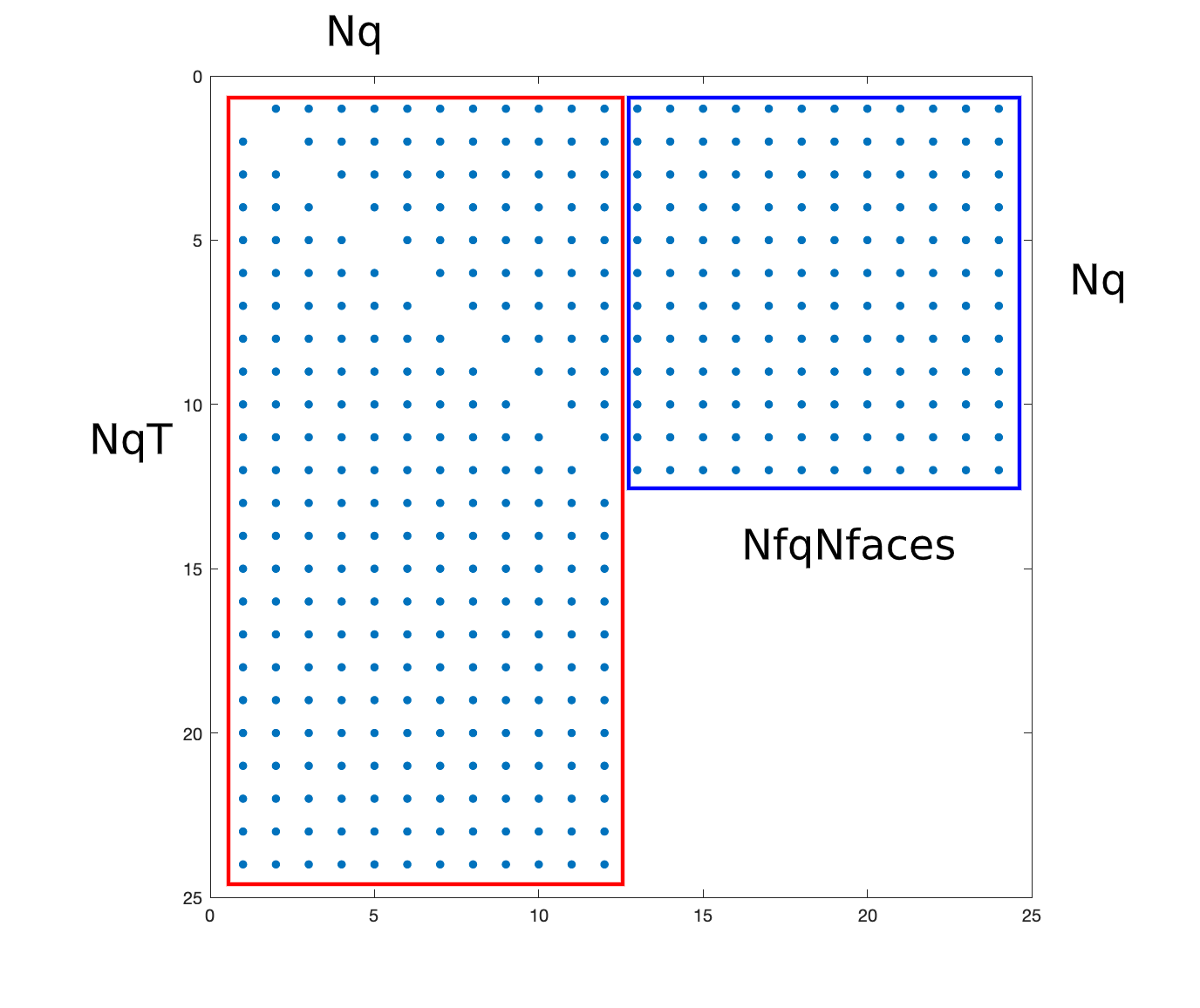}
    \caption{Spy diagram for $\bm{Q}^i_{h,skew} = \bm{Q}^i_h-(\bm{Q}^i_h)^T$}
    \label{fig:OCCA4}
    \end{figure}
\end{itemize}

\subsection{Performance comparison}
We conduct our GPU experiments on Google Cloud and provide the iteration time per step for K = 65536 elements. The results are presented in table \ref{table:GPUperformance}. We notice that the percentages of improvement cluster around 50\% for most of the experiments except for $N=1$. The ``baseline code'' includes optimization Step 1, but not Steps 2-4. 
\vspace{1.5em}
\begin{table}[t]
\begin{center}
\begin{tabular}{ |c|c|c|c|} 
 \hline
  N & Baseline & Optimized & Improvement percentage\\ 
 \hline
 1 & 0.0368s & 0.0221s & 39.94\%\\ 
 \hline
 2 & 0.0607s & 0.0316s & 47.94\%\\ 
 \hline
 3 & 0.0947s & 0.0515s & 54.38\%\\ 
 \hline
 4 & 0.1467s & 0.0627s & 57.26\%\\ 
 \hline
 5 & 0.3202s & 0.1585s & 50.50\%\\ 
 \hline
\end{tabular}
\end{center}
\caption{GPU runtime comparison between baseline and optimized version of implementation}
\label{table:GPUperformance}
\end{table}

\subsection{Performance comparison between ESDG and traditional DG}
The main difference between an entropy stable DG scheme and a traditional DG scheme is the volume kernels, which computes the DG approximation to a flux derivative. An entropy stable DG scheme computes a Hadamard product between two matrices, while a traditional DG method computes a regular matrix vector product. In order to explore how the GPU affects the computational efficiency of each approach, we compare both traditional and entropy stable volume kernels on both the CPU and GPU. The traditional volume kernel computes a dense matrix-vector product, while the entropy stable volume kernel performs flux differencing. Instead of actual SBP differentiation matrices, the kernels use square matrices with randomly generated entries in order to test computational performance for a wide range of matrix sizes. We note that these kernels most closely resemble volume kernels for traditional SBP schemes. The cost of the volume kernel for hybridized SBP operators will be higher due to the larger size of the matrices involved.

We define the run times to finish the entropy stable DG volume kernel and traditional DG volume kernel as $t_{ESDG}$ and $t_{DG}$, respectively. We denote the cost ratio $R_{CPU}$ as the ratio of $t_{ESDG}$ to $t_{DG}$ on the CPU, and define $R_{GPU}$ to be the same ratio for the GPU run time. We plot $R_{CPU}$ and $R_{GPU}$ for various matrix size in Figure \ref{fig:flux_GPU_serial_comp}.

\begin{figure}
\centering
\includegraphics[width=0.75\textwidth]{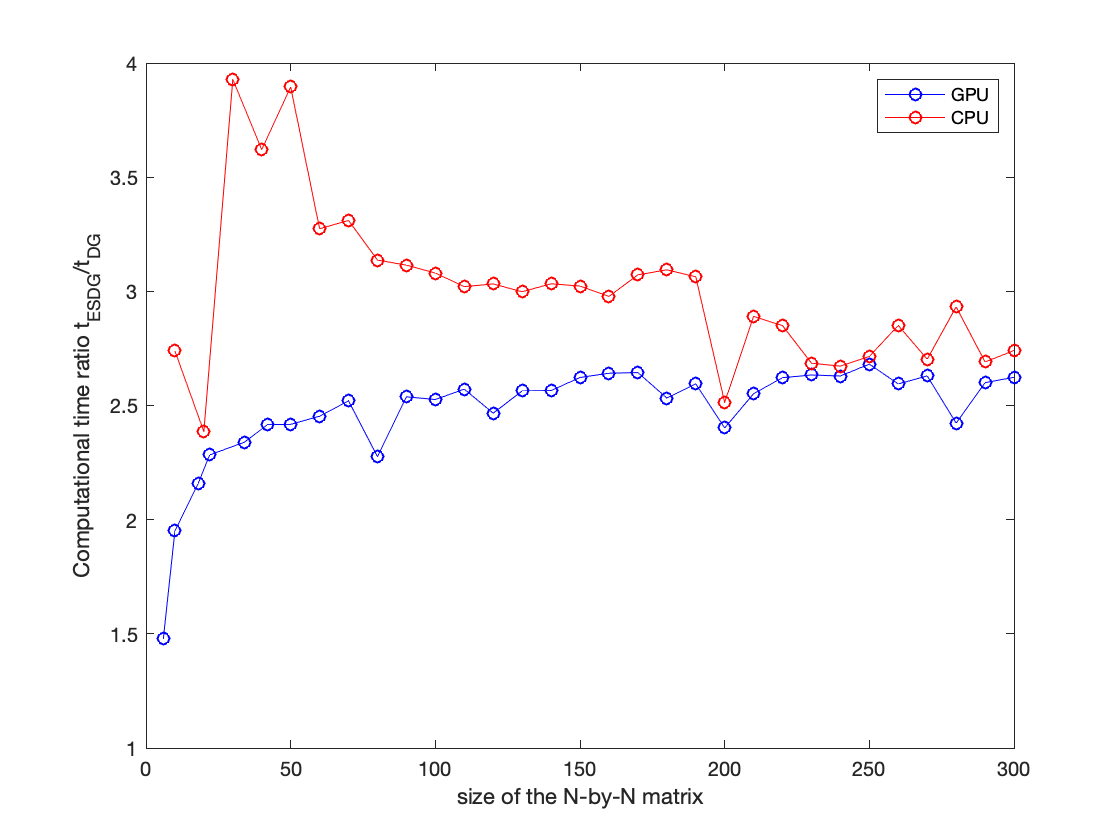}
\caption{Ratio of runtime between an entropy stable and traditional DG volume kernel.}
\label{fig:flux_GPU_serial_comp}
\end{figure}

We observe that the cost ratio is always lower for the GPU, but not very significantly. This implies that GPU implementation reduces the gap between ESDG and traditional DG computational times only slightly. For low to moderate degree $N$, GPU implementations perform better for the ESDG kernel relative to the CPU implementation. To interpret Figure \ref{fig:flux_GPU_serial_comp}, we point out that a $50 \times 50$ matrix size usually translates into a set of SBP nodes for a degree 5 or 6 polynomial approximation in 2D. The ratio drops again when the matrix size reaches around $200 \times 200$, which corresponds to a polynomial approximation of degree $N\approx 12$. 

The number of FLOPS for the computation of the term $(\bm{Q}_{h,skew}^i\circ \bm{F}^i)\bm{1}$ in our our triangular ESDG kernels is asymptotically the same as the number of operations needed to compute the matrix-vector product $\bm{Q}^i\bm{f}^i$ for the conventional DG method. Both involve operations involving $N_q^2$ matrix entries; however, we have $N_q^2$ flux evaluations in ESDG but only $N_q$ flux evaluations for the conventional DG method. In our ESDG implementation, we calculate entries of the flux matrix $\bm{F}_{ij}$ for each non-zero entry of $\bm{Q}^i_{h,skew}$.

Since the number of operations is $O(N_q^2)$ for both regular and ESDG, we expect the ratio of runtimes to eventually converge to a constant value. The ratio of 2.5 suggests that the cost of evaluating the flux for each non-zero matrix entry is 2.5x more expensive than computing contributions from a single matrix entry to a matrix-vector product for the conventional DG method.

This behavior is very different from the results shown in Figure 7 in \cite{wintermeyer2018entropy}, where the author found that routines of the traditional kernel and ESDG volume kernel were about the same for degree $N = 1,...,7$ on the GPU. The ESDG GPU implementation on quadrilateral meshes is memory-bound due to the fact that quadrilateral SBP operators are Kronecker products of 1D operators. Because of this, computational kernels can apply SBP operators in a more efficient manner using small $(N+1)\times (N+1)$ matrices corresponding to one-dimensional discretizations. The memory-bound nature of this kernel implies that memory transfers are the bottleneck. Thus, increasing the number of arithmetic operations does not increase the overall runtime until the cost of these operations exceeds the cost of memory traffic. As the polynomial order $N$ increases, the additional computation introduced in the ESDG volume kernel eventually increases the cost beyond that of the volume kernel for traditional DG. For $N$ sufficiently large, we expect the ratio of runtimes between ESDG and traditional DG on quadrilateral meshes to also approach a fixed value. However, it is not immediately clear that the ratio should approach 2.5, as this may depend on the hardware used and specific implementations.
On triangular meshes, our results show that the ESDG volume kernel is at least 1.5 times slower than the traditional DG volume kernel, and usually about 2.5 times slower for higher polynomial degrees. This difference arises due to the use the triangular meshes instead of quadrilateral meshes. Traditional volume kernels on triangles are not bandwidth bound \cite{chan2017bbdg}, and performing additional arithmetic operations results in a more significant increase in runtime. 
\end{subequations}
\section{Conclusions}
\label{sec:Conclusion}
\hspace{1em}In this work we present and compare two high-order entropy conserving and entropy stable schemes for the two dimensional shallow water equations on general curved triangular meshes. We construct well-balanced and high order entropy stable DG schemes for nonlinear conservation laws using hybridized SBP operators.  The resulting schemes satisfy a discrete conservation or dissipation of entropy. We compared entropy stable DG methods with DG-SBP methods based on traditional SBP operators, and compared the computational cost of entropy stable DG methods with traditional DG methods for both GPUs and CPUs.


\section*{Acknowledgement}
Authors Xinhui Wu and Jesse Chan gratefully acknowledge support from the National Science Foundation under awards DMS-1719818 and DMS-1712639. Author Ethan Kubatko gratefully acknowledges support from the National Science Foundation Grants ICER-1854991 and EAR-1520870. The authors also thank Benjamin Yeager for helping with the creation of the SBP quadrature rules.

\bibliographystyle{unsrt}
\bibliography{main.bib}

\end{document}